\documentclass{amsart}
\newtheorem{thm}{Theorem}[section]

\newtheorem{defn}[thm]{Definition}
\newtheorem{lemma}[thm]{Lemma}

\newtheorem{remark}[thm]{Remark}

\usepackage{amsmath}
\usepackage{amsxtra}

\newcommand{\id}{{\rm id}}
\usepackage{amscd}
\usepackage{amsthm}
\usepackage{amsfonts}
\usepackage{amssymb}
\usepackage{eucal}

\newcommand{\g}{{\mathfrak{g}}}
\renewcommand{\a}{{\mathfrak{a}}}
\newcommand{\h}{{\mathfrak{h}}}

\newcommand{\gl}{{\mathfrak{gl}}}
\newcommand{\n}{{\mathfrak{n}}}
\renewcommand{\sl}{{\mathfrak{sl}}}
\newcommand{\slr}{{\sl_{r+1}}}
\newcommand{\hg}{{\widehat{\g}}}

\newcommand{\ch}{{\rm ch}}
\newcommand{\G}{{\mathcal G_{\lambda,\boldsymbol\mu}(\boldsymbol\zeta)}}
\renewcommand{\a}{{\mathfrak a}}

\newcommand{\C}{{\mathbb C}}

\newcommand{\Z}{{\mathbb Z}}
\newcommand{\N}{{\mathbb N}}

\newcommand{\qbin}[2]
{{
\left[
\begin{matrix}{\displaystyle #1}\\
{\displaystyle #2}\end{matrix}
\right]
}}

\newcommand{\cF}{{\mathcal F}}
\newcommand{\cH}{{\mathcal H}}
\newcommand{\cK}{{\mathcal K}}

\newcommand{\cS}{{\mathcal S}}

\newcommand{\bN}{{\mathbf N}}

\newcommand{\bV}{{\mathbf V}}
\newcommand{\bW}{{\mathbf W}}

\newcommand{\bl}{{\mathbf l}}
\newcommand{\bm}{{\mathbf m}}
\newcommand{\bn}{{\mathbf n}}

\newcommand{\bx}{{\mathbf x}}
\newcommand{\by}{{\mathbf y}}

\newcommand{\al}{{\alpha}}
\newcommand{\bfalpha}{{\boldsymbol \alpha}}

\newcommand{\bomega}{{\boldsymbol{\omega}}}
\newcommand{\bal}{{\boldsymbol{\alpha}}}

\newcommand{\KK}{{\mathbb K}}

\newcommand{\NN}{{\mathbb N}}

\newcommand{\set}[2]{\{#1\ |\ #2\}}
\newcommand{\half}{{\frac{1}{2}}}

\newcommand{\bmu}{{\boldsymbol{\mu}}}
\newcommand{\bnu}{{\boldsymbol{\nu}}}

\newcommand{\olambda}{\overline{\lambda}}
\newcommand{\omu}{\overline{\mu}}

\newcommand{\xa}{x^{(\alpha)}}
\newcommand{\xb}{x^{(\beta)}}
\newcommand{\xaa}{x^{(\alpha+1)}}
\newcommand{\xap}{x^{(\alpha_p)}}
\newcommand{\ya}{y^{(\alpha)}}
\newcommand{\yb}{y^{(\beta)}}
\newcommand{\yaa}{y^{(\alpha+1)}}
\newcommand{\yap}{y^{(\alpha_p)}}
\newcommand{\ma}{m^{(\alpha)}}
\newcommand{\na}{n^{(\alpha)}}
\newcommand{\mua}{{\mu^{(\alpha)}}}
\newcommand{\nua}{{\nu^{(\alpha)}}}

\newcommand{\vbm}{{\overset{\to}{\bm}}}
\newcommand{\vbn}{{\overset{\to}{\bn}}}

\DeclareMathOperator{\mmod}{mod}
\DeclareMathOperator{\Gr}{Gr}
\DeclareMathOperator{\sym}{Sym}

\usepackage{graphicx}

\newcommand{\beq}{\begin{equation}}
\newcommand{\eeq}{\end{equation}}

\numberwithin{equation}{section}

\begin{document}
\title[Fermionic characters of $\widehat{\sl}_{r+1}$] {Fermionic characters
  and arbitrary highest-weight integrable $\widehat{\sl}_{r+1}$-modules}

\author{Eddy Ardonne, Rinat Kedem, Michael Stone}
\address{EA: Department of Physics, University of
    Illinois, 1110 W. Green St.,  Urbana, IL 61801. ardonne@uiuc.edu}
\address{RK: Department of 
    Mathematics, University of Illinois, 1409 W. Green Street, Urbana,
    IL 61801. rinat@uiuc.edu}
\address{MS: Department of
    Physics, University 
      of Illinois, 1110 W. Green St., Urbana, IL 61801.
\mbox{m-stone5@uiuc.edu}}
\begin{abstract}
  This paper contains the generalization of the Feigin-Stoyanovsky
  construction to all integrable $\widehat{\sl}_{r+1}$-modules.  We
  give formulas for the $q$-characters of any highest-weight
  integrable module of $\widehat{\sl}_{r+1}$ as a linear combination
  of the fermionic $q$-characters of the fusion products of a special
  set of integrable modules.  The coefficients in the sum are the
  entries of the inverse matrix of generalized Kostka polynomials in
  $q^{-1}$. We prove the conjecture of Feigin and Loktev
  regarding the $q$-multiplicities of irreducible modules in the
  graded tensor product of rectangular highest weight-modules in the
  case of $\sl_{r+1}$.  We also give the fermionic formulas for the
  $q$-characters of the (non-level-restricted) fusion products of
  rectangular highest-weight integrable $\widehat{\sl}_{r+1}$-modules.
\end{abstract}
\date{\today}

\maketitle
\section{Introduction}
Fermionic formul\ae\ for characters of highest-weight modules of affine
algebras or vertex algebras first appeared in a purely algebraic
context \cite{LP85}.  They were later shown \cite{KM,KKMM} to be
related to the partition functions of certain statistical mechanical
systems at their critical points.  These character formul\ae\ have
desirable combinatorial properties, such as the manifest positivity of
the coefficients that represent weight-space multiplicities. They also
have a physical significance because they reflect the quasi-particle
content of the statistical mechanical system. Consequently, algebraic
constructions of bases for representations which reveal this
combinatorial structure are important, and have been studied using
several methods in the past dozen years.

One such method is that of Feigin and Stoyanovski\u{\i} \cite{FS}. 
These authors used a theorem of Primc
\cite{Pr} to give an interesting construction of the vacuum integrable
modules of the affine algebra $\widehat{\g}$ associated to any simple
Lie algebra $\g$. Their construction relies on the loop generators of
the affine algebra.  Physical systems associated with such integrable
$\widehat{\g}$-modules are generalizations of the Heisenberg spin
chain in statistical mechanics, or the WZW model in conformal field
theory.

The formul\ae\ of Feigin-Stoyanovski\u{\i} \cite{FS} have an
attractive interpretation in terms of (a bosonic version of)
non-abelian quantum Hall states \cite{MR,AKS}.  In these states there
are $r$ ``types'' of particles that obey a generalized exclusion
principle: the wave function vanishes if any $k+1$ particles occupy
the same state. Here $r$ is the rank of the algebra and $k$ is the
level of the integrable $\widehat{\g}$-module. In the presence of
quasi-particle excitations, the wave functions can also vanish if
fewer than $k+1$ particles occupy the same state. The statistics of
the quasi-particles is `dual' to the statistics of the fundamental
particles \cite{ABS}.

The original construction of Feigin-Stoyanovski\u{\i} can be used to
compute \cite{FS} characters of vacuum (with highest weight
$k\Lambda_0$) representations of affine algebras. Later, Georgiev
\cite{Ge1,Ge2} generalized it to some modules in the ADE series,
with particularly simple highest weights, of the form
$l\omega_j+k\Lambda_0$, corresponding to special rectangular Young
diagrams.  (Here $\omega_j$ are certain fundamental $\g$-weights, and
$l\in\Z_{\geq 0}$.)  

In general, no fermionic formul\ae\ are available for arbitrary
highest-weight, integrable $\widehat{\g}$-modules. In this paper, we
resolve this problem for the case of $\sl_{r+1}$.

We explain, in terms of the functional realization of Feigin and
Stoyanovsk\u{i}, why such `rectangular highest weight' modules are
very special, and why there is no direct fermionic construction for
other modules.  However, we prove that it is possible to compute the
character of any module as a finite sum of fermionic characters of the
`rectangular' highest-weight modules. The coefficients in this sum are
the entries of the inverse matrix of generalized Kostka polynomials.
These coefficients are, however, not manifestly positive (or even of
positive degree).

In our construction we are naturally led to the graded tensor product
of Feigin and Loktev \cite{FL} of finite-dimensional $\g$-modules.  In
the case of irreducible $\sl_{r+1}$-modules with highest weights of the
form $l \omega_j$ (where $\omega_j$ is any fundamental weight), we
compute the explicit fermionic form of the graded multiplicities
of irreducible modules in the the Feigin-Loktev tensor product, thus
proving two of the conjectures of \cite{FL}: That the graded tensor product in
this case is independent of the evaluation parameters, and that it is
related to the generalized Kostka polynomials of \cite{SW,KS}.

The plan of the paper is as follows. In Section 2 we give the basic
definitions of the algebra and its modules. In Sections 3 and 4, we
supply the details of the generalized construction of \cite{FS} for
integrable modules of $\widehat{\sl}_{r+1}$, with highest weights
corresponding to rectangular Young diagrams. In Section 5, we explain
a similar calculation of graded characters of conformal blocks or
coinvariants (the fusion product of \cite{FL}), which turn out to be
related to the generalized Kostka polynomials of \cite{SW,KS}. We then
use this calculation in Section 6 to compute the characters of
arbitrary highest-weight representations. See theorem \ref{chvl} for
the main result.

Although, for the sake of clarity, we concentrate in this paper on the
case of $\widehat\g=\widehat\sl_{r+1}$, the generalization to affine
algebras associated with other simple Lie algebras is possible, but in
that case one should replace the notion of integrable
$\widehat{\g}$-modules with irreducible $\g$-modules as their top
component with those which have (the degeneration to the classical
case of) Kirillov-Reshetikhin modules as their top component. We will
give this construction in a future publication.

{\em Acknowledgements:}
The work of E.A. is supported by NSF grants numbers DMR-04-42537 and
DMR-01-32990; that of M.S. by NSF grant DMR-01-32990.
R.K. would like to thank B. Feigin and S. Loktev
for many useful discussions.

\section{Notation}
\subsection{Current generators of affine algebras}
Let $\g=\sl_{r+1}$ and let $\Pi=\{\alpha_i \ |\ i=1,\ldots,r\}$ denote its
simple roots, and $\{\omega_i \ |\ i=1,\ldots,r\}$
the fundamental weights. Let $\{e_{\alpha_i}=e_{i}\ |\ i=1,\ldots,r\}$
denote the corresponding generators of
$\n_+$, and $\{f_{\alpha_i}=f_{i}\ |\ i=1,\ldots,r\}$ those of $\n_-$. 
We have the Cartan
  decomposition $\sl_{r+1} \simeq \n_+ \oplus \h\oplus \n_-$, where
$\h$ is the Cartan subalgebra.
  
  Irreducible, finite-dimensional highest-weight $\g$-modules
  $\pi_\lambda$ are parametrized by weights $\lambda\in P^+$, that is,
  $\lambda = l_1 \omega_1 + \cdots + l_{r}\omega_{r}$  with
  $l_i\in\Z_{\geq 0}$. The subset of $P^+$ consisting of weights
  $\lambda$ such that 
  $\sum_{i=1}^{r} l_i \leq k$ is called the set of level-$k$
  restricted weights, $P^+_k$.

The affine Lie algebra associated with $\g$ is $\widehat{\g}$, where
$$
\widehat{\g}\simeq \g\otimes \C[t,t^{-1}]\oplus \C c\oplus \C d,
$$
where $c$ is central and
\begin{equation} \label{degreedef}
[d,x\otimes t^n]=-n x\otimes t^n \ .
\end{equation}
We denote the current generators by $x[n]\overset{\rm def}{=}
x\otimes t^n$, $x\in \sl_{r+1}$. Let $\langle x,y\rangle$ be the
symmetric bilinear form on $\slr$. Then the relations between the
currents are
$$
[x\otimes f(t),y\otimes g(t)]_{\widehat{\g}} = [x,y]_{\g} f(t) g(t) +
c\ \langle x,y\rangle \oint_{t=0} f'(t) g(t) dt,
$$
where $[\cdot,\cdot]_{\g}$ is the corresponding commutator in $\g$.

The Cartan decomposition is
$\widehat{\g}\simeq \widehat{\n}_+ \oplus \widehat{\h} \oplus
\widehat{\n}_-$ with
$\widehat{\n}_\pm = \n_\pm \oplus (\sl_{r+1}\otimes
t^{\pm1}\C[t^{\pm1}])$
and $\widehat{\h}= \h \oplus \C c \oplus \C d$.
The algebra $\widehat{\g}'$ is the algebra obtained by dropping the
generator $d$.

We will frequently use generating functions for current generators of
the affine algebra, which we define by
\begin{equation}\label{current}
x(z) = \sum_{n\in \Z} x[n] z^{-n-1},\ x\in\sl_{r+1}.
\end{equation}
Note that the convention for the current generators in \eqref{current}
is different from that used by \cite{FS,FJKLM}.

\subsection{Affine algebra modules}

On any irreducible $\widehat{\sl}_{r+1}$-module, $c$ acts by a constant
$k$ called the {\em level} of the representation. A cyclic highest-weight
$\hg$-module with highest weight $\Lambda=\lambda + k\Lambda_0 +
m\delta$ is a cyclic module generated by the action of $\hg$ on a
highest-weight vector $v_\lambda$, such that 
\begin{eqnarray}
&& \widehat{\n}_+ v_\lambda = 0, \label{highest}\\
&& h v_\lambda = \lambda(h) v_\lambda, \ {\rm for}\  h\in \h\subset\g, \quad
 c v_\lambda = k v_\lambda, \quad
 d v_\lambda = m v_\lambda.\label{weight}
\end{eqnarray}
The universal such module is the Verma module $M(\Lambda)\simeq
U(\widehat{\n}_-)$. If $k\in 
\N$ and $\lambda\in P_k^+$, the quotient of the Verma module by its
maximal submodule is an irreducible, highest-weight integrable
$\hg$-module, which we denote by $V_\lambda$ (we assume $k$ is fixed
in this notation).  The structure of the cyclic module generated by a
highest-weight vector $v_\lambda$ is independent of $m$, so it is
generally convenient to set $m=0$.

\begin{defn}
  Let $M$ be an irreducible cyclic highest-weight module with highest
  weight $\Lambda = \lambda + k\Lambda_0$, generated by
  the highest-weight vector $v_\lambda$. The subspace generated by the
  action of the subalgebra $\g\otimes 1\simeq \g$ on $v_\lambda$ is
  called the {\em top component} of $M$. It is isomorphic as a
  $\g$-module to $\pi_\lambda$.
\end{defn}

The irreducible, finite-dimensional $\g$-module $\pi_\lambda$ is
characterized as the quotient of the Verma module of $\g$ by the left
ideal in $\g$ generated by $f_i^{l_i+1}$.  Similarly, the integrable
module $V_\lambda$ is the quotient of the Verma module $M(\Lambda)$ of
$\widehat{\g}$ by the left ideal in $\hg$ generated by
$f_i[0]^{l_i+1}$, plus one additional generator,
$e_\theta[-1]^{k-\theta(\lambda)+1}$ where $\theta=\alpha_1 + \cdots +
\alpha_r$.

A characterization of the maximal proper submodule $M'(\Lambda)$ of
$M(\Lambda)$ in the case of integrable modules was given in \cite{Pr}
in terms of the algebra of current generators.

Note that on any highest-weight module, the current
(\ref{current}) acts as a Laurent series in $z$. Therefore, products
of currents make sense when acting on a highest-weight module, and one
can consider the associative algebra of currents.  Formally, the
coefficients of $z^n$ in products of currents of the form $x(z) y(z)$
exist only in a completion $\overline{U}$ of $U(\widehat{\g})$.

\begin{thm}\label{primc}\cite{Pr}
  Let $M(\Lambda)$ be a Verma module with highest weight $\Lambda =
  \lambda + k \Lambda_0$, with $\lambda\in P^+_k$ and $k\in \N$. Denote
  its maximal proper submodule by $M'(\Lambda)$, such that
  $V_\lambda \simeq M(\Lambda)/M'(\Lambda)$. Let $R$ be the subspace in
  $\overline{U}$ generated by the adjoint action of $U(\sl_{r+1})$ on the
  coefficients of $e_\theta(z)^{k+1}$. Then $M'(\Lambda) = R
  M(\Lambda)$.
\end{thm}
Again, the elements in $R$ act as well-defined elements of
$U(\widehat{\g})$ on $M(\Lambda)$. We call the set of currents which
result from the adjoint action of $\sl_{r+1}$ on the current
$e_\theta(z)^{k+1}$ the {\em integrability conditions}. For example,
for any root $\alpha$, the coefficients of $e_\alpha(z)^{k+1}$ are in
$R$.

\section{The semi-infinite construction of Feigin and
  Stoyanovski\u{\i}}\label{feig-sto}
Theorem \ref{primc} was used by Feigin and Stoyanovski\u{\i} \cite{FS} to
give a construction of the integrable modules in the case where
$\Lambda=k\Lambda_0$.  The construction naturally gives rise to
fermionic formul\ae\ for the characters of integrable modules.  We will
explain the details of the construction of \cite{FS} below.

\subsection{Principal subspaces}
\newcommand{\wn}{{\widetilde{\n}_-}}
For arbitrary integrable highest weight
$\Lambda=\lambda+k \Lambda_0$, let $v_\lambda$ be
the highest-weight vector of $V_{\lambda}$. Consider the subalgebra
$$\widetilde{\n}_-\overset{\rm def}{=}\n_-\otimes \C[t,t^{-1}]$$
acting on $v_\lambda$.  

\begin{defn} Define the principal subspace $W_\lambda =
  W_\lambda^{(0)} = U(\wn) v_\lambda\subset V_\lambda$. Similarly,
  define the principal subspaces $W_\lambda^{(N)}= U(\wn) T_{N}
  v_\lambda$, where where $T_{N}=t_{\alpha(N)}$ is the affine Weyl
  translation corresponding to the root $\alpha(N)=\sum_i N_i \alpha_i$
  (in the notation of \cite{Kac} (6.5.2)), where $N_i$ are positive
  integers such that $(C_r\bN)_\al = 2N$ for all $\al$, and $C_r$ is
  the Cartan matrix of $\sl_{r+1}$.
\end{defn}

\begin{lemma}
This choice of $\alpha(N)$ gives a sequence of inclusions
\begin{equation}\label{incl}
W_\lambda^{(0)} \subset W_\lambda^{(1)} \subset \cdots \subset
W_\lambda^{(N)} \subset \cdots,
\end{equation}
such that the inductive limit of the sequence (\ref{incl}) as
  $N\to\infty$ is the integrable module $V_\lambda$.
\end{lemma}

The inclusions follow from the fact that $v_\lambda\in
W^{(N)}_\lambda$. The fact that the inductive limit indeed gives the
full module is not obvious (see \cite{Primc,FJLMM}) but follows from
the fact that the module is integrable.

In fact, this Theorem was proven in \cite{FS} for the following cases:
$\widehat{\sl}_2$ for arbitrary highest weight, and $\widehat{\sl}_3$
with $\Lambda = k \Lambda_0$. This was done by computing the
characters in the limit $N\to\infty$, and comparing them with the
known character formul\ae\ for $V_\lambda$ of \cite{LP85}.

In \cite{Ge1,Ge2}, certain combinatorial proofs were provided using
ideas related to those of \cite{FS} (with differently defined
principal subspaces) for rectangular highest weights, for all simply
laced algebras. The principal subspaces of that paper are different
from those used here, as \cite{Ge1} uses what amounts to a different
subalgebra to generate the subspace.

In this paper, we will continue this program by giving the character
formul\ae\ for arbitrary highest-weight modules of $\sl_{r+1}$. It
turns out that the methods of \cite{FS} are not sufficient for the
case of non-rectangular representations, and instead we must resort to
computing the characters of certain fusion products of
representations, and decomposing them in terms of irreducible
modules. The result is a formula which is a sum of fermionic formulas
of the form found in \cite{FS,Ge1,Ge2}, where the coefficients in the
sum are elements of $\Z[q^{-1}]$.

\subsection{Relations in the principal subspace}
Let us characterize the ideal $I_\lambda$, where $W_\lambda \simeq
U(\wn)/I_\lambda$. 
Using a PBW-type argument, it is easy to see that $W_\lambda =
U(\n_-\otimes \C[t^{-1}] ) v_\lambda$, because the highest-weight
vector $v_\lambda$ is annihilated by $\n_-\otimes t \C[t]$.  Thus,
$I_\lambda$ includes the left ideal generated by $\{ f_\alpha[n]\ |\ n>0,
\alpha\in \Pi\}$.

The ideal contains the two-sided ideal generated by relations in the
Lie algebra. In terms of generating functions, these relations are
\begin{eqnarray}
[f_{\alpha_i}(z),f_{\alpha_j}(w)] &=& \left\{\begin{array}{ll} 0, &
    |i-j|\neq 1 \\
w^{-1}\delta(w/z) f_{\alpha_i+\alpha_j}(z), & |i-j|=1
\end{array}\right. \label{comm} 
\\
\left[f_{\alpha_i}(z),{[f_{\alpha_i}(w),f_{\alpha_{i\pm1}}(u)]}\right]
&=& 0 \label{serre},
\end{eqnarray}
where $\delta(z) = \sum_{n\in\Z} z^n$.
These two relations together mean that matrix elements involving the product
$f_{i}(z)f_{{i\pm1}}(w)$ have a simple pole whenever
$z= w$, and that the residue of this pole commutes with
$f_{i}(u)$. 

The integrability condition 
\begin{equation}\label{integ}
f_{i}(z)^{k+1}v=0,\quad v\in V_\lambda,\ 1\leq i \leq r,
\end{equation} for any
implies that $I_\lambda$ contains the
two-sided ideal generated by the coefficients of $z^n$ of
$f_{i}(z)^{k+1}$ (in the appropriate completion of the
universal enveloping algebra).

Finally there are the relations which follow from the integrability of
the top component $\pi_\lambda$ of $V_\lambda$, which is a subspace of
$W_\lambda$ also. Therefore, $I_\lambda$ contains the left ideal
generated by $f_{i}[0]^{l_i+1}$.  The integrability condition
involving $e_\theta[-1]$ does not play a role, because it is not an
element of $U(\wn)$.

\subsection{Construction of the dual space}
In order to compute the characters of the principal subspace
$W_\lambda$, we describe its dual space. This will enable us to
calculate the character for sufficiently simple $\lambda$. The dual
space is spanned by the coefficients of monomials of the form
$x_1^{n_1}\cdots x_m^{n_m}$ of matrix elements in the set
$$
\mathcal F_{\lambda} = \left\{\langle w |
f_{{i_1}}(x_1)\cdots f_{{i_m}}(x_m)| v_\lambda\rangle \ |\
\quad w\in V_\lambda^*,\ m\geq 0,\ 1\leq i_a\leq r \right\},
$$
where $V^*_\lambda$ is the restricted dual module. 
Given an ordering of the generators, the function above is defined in
the region $|x_i|>|x_{i+1}|$, and therefore the coefficient of
$x_1^{n_1}\cdots x_m^{n_m}$ for given integers $n_j$ is given by the
expansion in this regime. Below, we shall refer to the function space
$\mathcal F_\lambda$ itself as the dual space, and specify an
appropriate pairing. This space can be characterized by its pole structure
and vanishing conditions.

\subsubsection{The dual space to $U(\wn)$}
Let us first consider the larger function space $\mathcal G$, dual to
the universal enveloping algebra $U=U(\wn)$.  The algebra $U$ is
spanned by words in the letters $\{f_{\alpha_i}[n]\ |\  i=1,\ldots,r,\ n\in
\Z\}$, and it is $\h$ and $d$-graded. The graded component
$U[\mathbf m]_d$, where $\bm=(m^{(1)},\ldots,m^{(r)})^T$,
is spanned by the elements
$f_{i_1}[n_1]\cdots f_{i_m}[n_m]$, of $\h$-weight $\sum_\alpha
m^{(\alpha)}\alpha= \sum_j \alpha_{i_j}$ and $-\sum_i n_i = d$.

\newcommand{\cG}{{\mathcal G}}
The dual space to $U$ is also $\h$- and $d$-graded. Denote by $U[\bm]$
the $\h$-graded component, and by $\cG[\mathbf m]$ the dual to it.
This is a space of functions in the variables 
$$\mathbf x = \{x_i^{(\alpha)} \ |\
i=1,\ldots,m^{(\alpha)}, \alpha=1,\ldots,r\},$$
where $x_i^{(\alpha)}$ is the variable corresponding to a generator of
the form $f_\alpha(x_i^{(\alpha)})$.
We define the pairing $(\cdot,\cdot)$ between $U$ and $\mathcal G$
inductively, as follows:
\begin{eqnarray}
&& ( 1,1) = 1\nonumber \\
&& (g(\mathbf x),  M f_\alpha[n] ) =
\left(\frac{1}{2\pi i}\oint_{x^{(\alpha)}_1=0}
(x_1^{(\alpha)})^n 
  g(\mathbf x)dx_1^{(\alpha)},M \right) 
,\quad M\in U,\label{pairing}
\end{eqnarray}
where the contour of integration is taken counter-clockwise around the
point $x^{(\alpha)}_1=0$, in such a way that all other points are excluded,
$|x_1^{(\alpha)}|<|x_j^{(\alpha')}|$.  Similarly,
\begin{eqnarray*}
&&(g(\mathbf x),   f_\alpha[n] M ) =
\left(\frac{1}{2\pi i}\oint_{x^{(\alpha)}_1=0}
(x_1^{(\alpha)})^n 
  g(\mathbf x)dx_1^{(\alpha)},M \right) ,
\end{eqnarray*}
the contour is taken clockwise.

The commutation relations between the currents are equivalent to the
operator product expansion (OPE)
$$
f_i(z) f_{i\pm 1}(w) = \frac{f_{\alpha_i+\alpha_{i\pm 1}}(w)}{z-w} +
\hbox{regular terms}
$$
where ``regular terms'' refers to terms which have no pole at
$z=w$, and the expansion of the denominator is taken in the region $|z|>|w|$.
Due to the OPE's, it
is clear that functions in $\mathcal G[\mathbf m]$ will have at most
a simple pole whenever $x_j^{(\alpha)} = x_k^{(\alpha\pm1)}$. Thus,
functions in $\mathcal G[\mathbf m]$ are rational functions of the
form 
\begin{equation}\label{rational}
g(\mathbf x) = \frac{g_1(\mathbf x)}{\prod_{i,j,\alpha}
  (x_i^{(\alpha)}-x_j^{(\alpha+1)})},
\end{equation}
where $g_1(\mathbf x)$ are polynomials in $(x_i^{(\alpha)})^{\pm1}$.

Again using the OPE's, we can construct the pairing between all
other elements of $U$ and $\mathcal G$. For example,
\begin{eqnarray*}
& &(g(\mathbf x), M f_{\alpha + \alpha\pm 1}[n] )
=\left( \frac{1}{2\pi i}\oint_{x_1^{(\alpha)}=0}
\left.  (x_1^{(\alpha)})^n(x_1^{(\alpha)}-x_1^{(\alpha\pm 
  1)}) g(\mathbf x)\right|_{x_1^{(\alpha)}=x_1^{(\alpha\pm1)}}
dx_1^{(\alpha)}, M\right) 
\end{eqnarray*}
where the contour excludes all other points, and
\begin{equation*}
\begin{split}
& (g(\mathbf x), M f_{\alpha +\cdots+ \alpha+h}[n] )=\\
& \left( \frac{1}{2\pi i}\oint_{x_1^{(\al)}=0}
 \left. (x_1^{(\al)})^n(x_1^{(\al)}-x_1^{(\al+1)})
\cdots(x_1^{(\al+h-1)}-x_1^{(\al+h)})
 g(\mathbf x)\right|_{x_1^{(\al)}=\cdots=x_1^{(\al+h)}}dx_1^{(\al)},
M\right).
\end{split}
\end{equation*}

The function $g_1(\mathbf x)$ is not completely arbitrary, due to
the Serre relation \eqref{serre}. The Serre relation implies that the function
\begin{equation*}
\left.(x_1^{(\al)} - x_1^{(\al+1)}) g(\mathbf
x)\right|_{x_1^{(\al)}=x_1^{(\al+1)}} 
\end{equation*}
has no poles at the points $x_j^{(\al+1)}=x_1^{(\al)}$ and
$x_j^{(\al)}=x_1^{(\al+1)}$, where $j>1$. This implies that the function
$g_1(\mathbf x)$ has the property that
\begin{equation}\label{serre1}
\left.g_1(\mathbf x) \right|_{x_i^{(\al)} = x_j^{(\al)} = x_k^{(\al\pm1)}} =
0.
\end{equation}

Finally, it is clear that since $[f_i(z),f_i(w)]=0$, $g_1(\mathbf
x)$ is symmetric under the exchange of variables
$x_i^{(\al)}\leftrightarrow x_j^{(\al)}$. In summary, we have
\begin{thm}\label{dualtoU}
The space of functions $\mathcal G[\mathbf m]$ dual to the graded component
$U[\mathbf m]$ of the universal enveloping algebra of
$\wn$, with the pairing defined inductively by
(\ref{pairing}), is the space of functions in the variables
$\{x_j^{(\al)}\}$ with $j=1,\ldots,m^{(\al)}$ and $\al=1,\ldots,r$, 
of the form (\ref{rational}), where $g_1(\mathbf x)$ is a polynomial
in $(x_j^{(\al)})^{\pm 1}$, symmetric under the exchange of variables
with the same superscript, and which vanishes whenever
$x_1^{(\al)}=x_2^{(\al)}=x_1^{(\al\pm1)}$.
\end{thm}

\subsubsection{Dual to the principal subspace $W_\lambda$}
Next, we consider the space $\mathcal F_\lambda[\mathbf m]$, which is
defined as the graded component of the space $\mathcal F_\lambda$, the
subset of matrix elements of $U[\mathbf m]$ in $\cF_\lambda$. The
space $\cF_\lambda[\bm]$ is
the dual space to $W_\lambda[\bm]$ (the weight subspace of $W_\lambda$
of $\h$-weight $\lambda-\bm^T\bal$) with the pairing defined as in
(\ref{pairing}), where $1\in U$ is replaced by $v_\lambda$.

The dual space $\mathcal F_{\lambda}[\mathbf m]$ is the subspace of
$\mathcal G[\mathbf m]$, which couples trivially via the pairing
(\ref{pairing}) to the ideal $I_\lambda\subset U$.  Apart from the
two-sided ideal coming from the relations in the algebra, which we
have already accounted for in constructing $\cG[\bm]$, the ideal
$I_\lambda$ contains the relations coming from the highest-weight
conditions (\ref{weight}), and from the integrability conditions
(\ref{integ}).

The integrability conditions mean that $U f_i(x)^{k+1} U\subset
I_\lambda$, which 
means that
\begin{equation}\label{simpleroot}
g_1(\mathbf x)|_{x_1^{(\al)} = \cdots = x_{k+1}^{(\al)}} = 0,
\end{equation}
for all $g(\mathbf x)\in \mathcal F_{\lambda}[\mathbf m]$ and for all
$\al$. 

The ideal $I_\lambda$ contains the left ideal generated by
$f_\al[n], n>0$ for any  $\al$.
We see from (\ref{pairing}) that for
functions in $\mathcal F_\lambda[\mathbf m]$, $g_1(\mathbf x)$ can
have at most a simple pole at $x_1^{(\al)}=0$. Let us define the
function $g_2(\mathbf x)$ by
\begin{equation}\label{gtwo}
g(\mathbf x) = \frac{g_2(\mathbf x)}
{\prod_{\al,i}(x_i^{(\al)})\prod_{\al,i,j}(
x_i^{(\al)}-x_j^{(\al+1)})} \ ,
\end{equation}
where  $g_2(\mathbf x)$ is a polynomial in $x_i^{(\al)}$ for all
$i,\al$.

In order to account for the relation $U f_\beta[n]\subset I_\lambda$ for
$\beta = \alpha_i + \cdots + \alpha_{i+h}$, where $n>0$, we need to
impose an additional restriction on $g_2(\mathbf x)$, because of the
prefactor $(x_1^{(\al)} x_1^{(\al+1)}\cdots x_1^{(\al+h)})^{-1}$ in
(\ref{gtwo}). The function $g_1(\mathbf x)$, after evaluation at the
point $u=x_1^{(\al)}=x_1^{(\al+1)}=\cdots=x_1^{(\al+h)}$, must be of degree
greater than or equal to $-1$ in the variable $u$ if it is to couple
trivially to $f_\beta[n]$ for $n>0$. Therefore, we see that
$g_2(\mathbf x)$ satisfies: 
\begin{equation}\label{compoundroot}
g_2(\mathbf
x)|_{x_1^{(\al)}=x_1^{(\al+1)}=\cdots=
x_1^{(\al+h)}=u}\quad \hbox{vanishes as $u^h$ as
$u\to 0$}.
\end{equation}

Finally we need to take into account the integrability conditions
for the top component: $U f_\beta[0]^{\lambda(\beta)+1}\subset
I_\lambda$ for each positive root $\beta$. For simple roots, this means that 
\begin{equation}\label{simple}
g_2(\mathbf x)|_{x_1^{(\al)} = \cdots = x_{l_\al+1}^{(\al)} = 0} = 0.
\end{equation}
When $\beta$ is not a simple root, then the relations are more
complicated, involving variables corresponding to different roots. 
These are sufficiently complicated that we do not know how to compute
the character of the space in this case.

However, at this point let us note that for the special
case of rectangular representations, the situation is much
simpler. The relation (\ref{compoundroot}) is automatically satisfied
for such representations. For suppose we consider the representation
with $l_\beta\neq 0$ for at most one index $\beta$. Then since
$Uf_\beta[0]\not\subset I_\lambda$, whereas $Uf_\al[0]\subset
I_\lambda$ for $\al\neq \beta$, we have that in this special case,
\begin{equation}\label{rectangular}
g_1(\mathbf x) = \prod_j (x_j^{(\beta)})^{-1} g_2(\mathbf x)
\end{equation}
where $g_2(\mathbf x)$ is a polynomial in all the variables,
satisfying (\ref{simple}) for the index $\beta$ only, as well as the
integrability conditions and the Serre relation.  The relation
(\ref{compoundroot}) is not an extra condition in this case.

Let us summarize the result for rectangular representations,
therefore.
\begin{thm}
  Let $\Lambda_\beta = l\omega_\beta + k \Lambda_0$ 
  for some $1\leq \beta\leq r$. Then the dual space of functions to the graded
  component of the principal subspace $W_{l\omega_\beta}[\mathbf m]$ is the
  space of rational functions of the form (\ref{rational}), where
  $g_1(\mathbf x)$ is a function of the form (\ref{rectangular}),
  where $g_2(\mathbf x)$ is a polynomial in the variables $x_i^{(\al)}$
  satisfying the Serre relation (\ref{serre1}), symmetric under the
  exchange of variables $x_i^{(\al)}\leftrightarrow x_j^{(\al)}$ for all $\al$,
  vanishing when $x_1^{(\beta)}=\cdots=x_{l+1}^{(\beta)} = 0$, or when
  any $k+1$ variables of the same superscript coincide, $x_1^{(\al)} = \cdots =
  x_{k+1}^{(\al)}$ for any $\al$.
\end{thm}
In the next section, we will show how to compute the character of this
space using a filtration on the space.

For non-rectangular representations there is no such simple
description of the space. The purpose of this paper is to explain how
to compute the character for non-rectangular representations as a
linear combination of characters of rectangular representations.

\subsection{Filtration of the dual space $\mathcal F_\lambda$}
\label{fildual}

In this subsection, we will assume that
$\Lambda=\Lambda_\beta = \lambda+ k \Lambda_0$,
$\lambda = \lambda_\beta=l\omega_\beta$ for some fixed $1\leq
\beta\leq r$. This 
corresponds to a Young diagram of rectangular form (with $l$ columns
and $\beta$ rows).

As explained above, the space $\cF_\lambda$ is $\h$-graded,
$\cF_\lambda = \bigoplus_\bm \cF_\lambda [\bm]$, where
$\mathcal F_\lambda [\bm]$
is a subspace of the space of rational functions in the variables
$\bx  = \set{\xa_i}{\alpha=1,\ldots,r \ ; \ i=1,\ldots,m^{(\alpha)}}$
of the form
\begin{equation}
\label{rationalf}
G(\bx) = \frac{g(\bx)}
{\prod_{i} (\xb_i) \prod_{\alpha=1}^{r-1}\prod_{j,k}
(\xa_j-\xaa_k)} \ ,
\end{equation}
where $g(\bx)$ is polynomial, symmetric under exchange of
variables with the same value of $\alpha$ (which we will refer
to as the {\em color} index),
$\xa_i \leftrightarrow \xa_{j}$. The index $\beta$ corresponds to
the to the fundamental weight $w_\beta$, where $\lambda_\beta= l w_\beta$.
In addition, 
$g(\bx)$ vanishes when any of the following conditions is met
\begin{align}
&\xa_1=\cdots=\xa_{k+1} \label{intc} \\ \label{serrec}
&\xa_1=\xa_2=x_1^{(\alpha\pm 1)} ,\\
\label{irrepc}
&\xb_1=\cdots=\xb_{l+1}=0 \ .
\end{align}

Our goal is to compute the character of this space, for which purpose
we will introduce a filtration and an associated graded space. We will
be able to compute the characters of the graded pieces easily.

To simplify the calculations below, let us define the closely related space
$\overline{\cF}_\lambda [\bm]$.  This space is a subspace the space of
all rational functions in the variables $\bx$, which are given by
\begin{equation}
\label{ratfb}
\overline{G}(\bx) = \frac{{g}(\bx)}
{\prod_{\alpha=1}^{r-1}\prod_{i,j} (\xa_i-\xaa_j)} \ ,
\end{equation}
where ${g}(\bx)$ is as in \eqref{rationalf}, so
$\overline{G}(\bx) = \prod_i \xb_i G(\bx)$.
In the following, we will fix $\bm$ and $l$, and 
study a filtration of this space
$\overline{\cF}_{\lambda}[\bm]$ (which we will refer to by
$\overline{\cF})$, which can be described as follows.

Let $\boldsymbol \mu = (\mu^{(1)},\ldots,\mu^{(r)})$ be a collection of
partitions, where each $\mu^{(\alpha)}$ is a partition of $m^{(\alpha)}$
and has $m^{(\alpha)}_a$ rows of length $a$.

\begin{figure}[ht]
\begin{center}
\includegraphics[width=5cm]{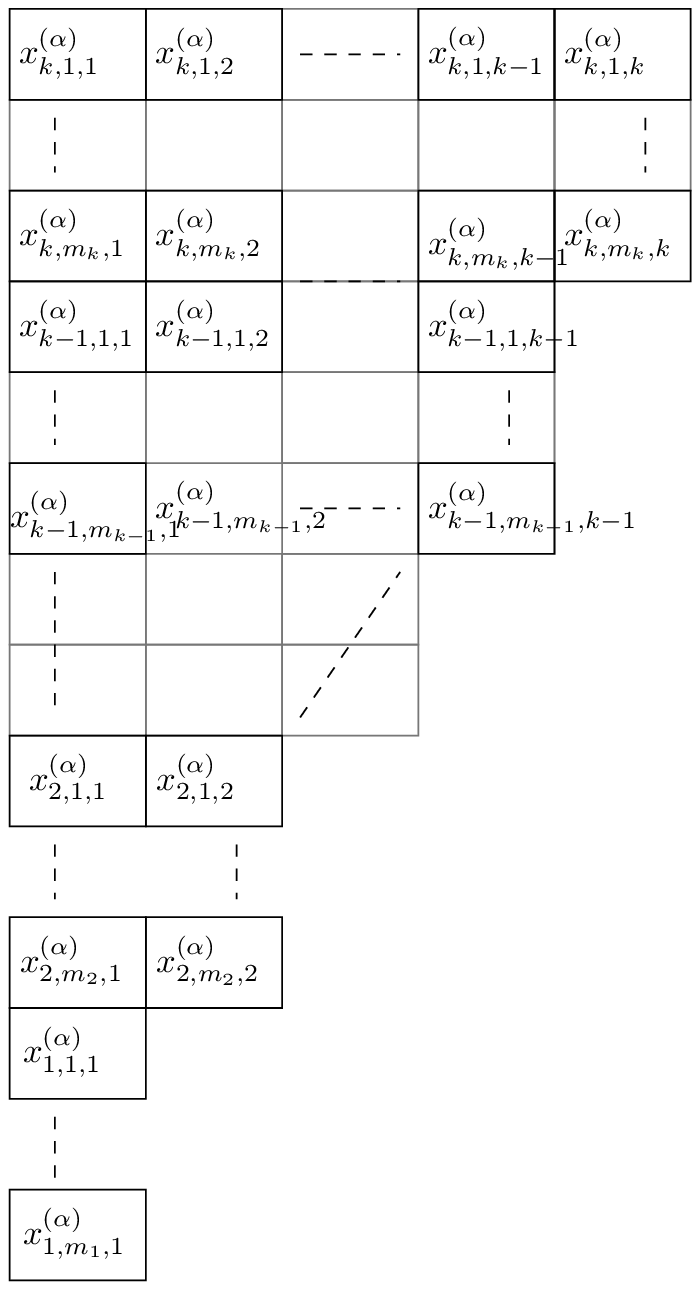} \qquad
\raisebox{5 cm}{$\longrightarrow$} \qquad
\includegraphics[width=5cm]{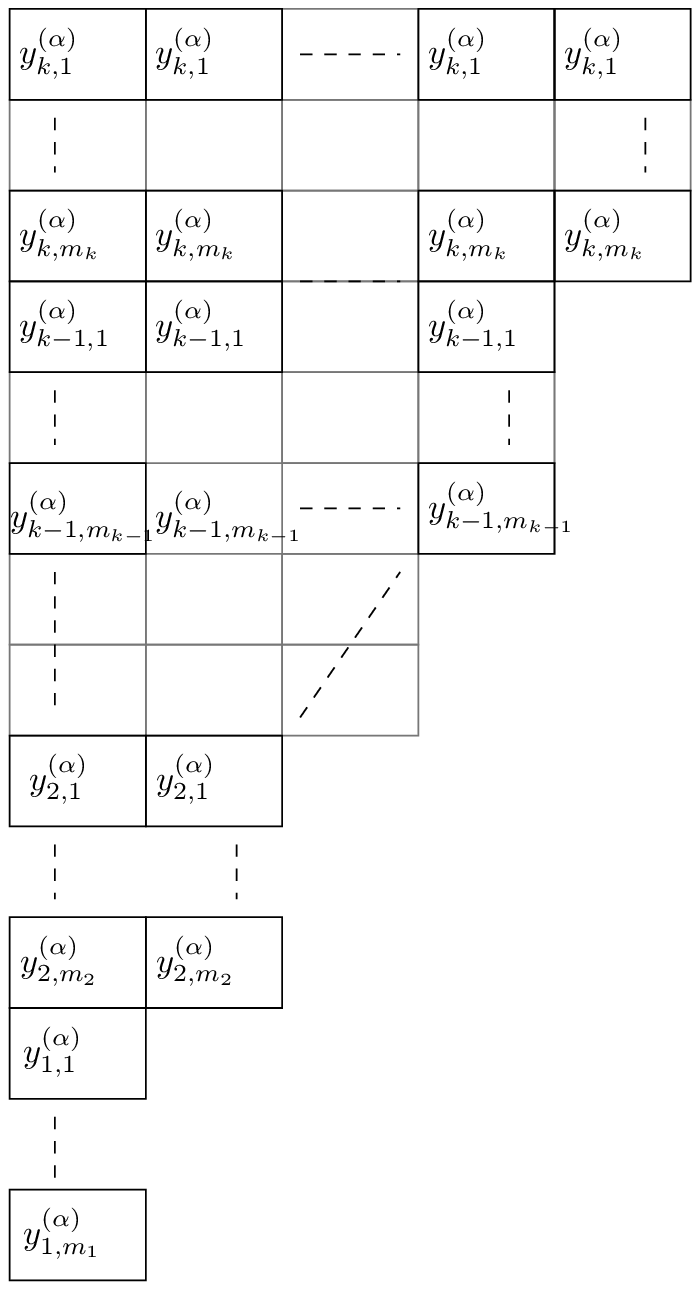}
\end{center}
\caption{The evaluation map for the variables $\xa$.
Note that we dropped the superscripts $(\alpha)$ in
$\ma_a$.}
\label{evaluationmap2}
\end{figure}

We can now rename the variables $\xa_i$ by associating 
each of them to a
box of the Young diagram associated with the partitions
$\mua$. As a result of this renaming, we have variables
$\xa_{a,i,j}$, which corresponds to the Young diagram of partition
$\mua$, namely to column $j$ of the $i^{th}$ row
(counted from top to bottom) of length $a$. See the left
part of figure \ref{evaluationmap2} for an explicit example.
In the proofs which follow, we will simplify this notation as much as
possible. Note that, due to the symmetry properties of $g(\bx)$, how
we rename the variables is irrelevant.  

Let $\cH$ be the space of rational functions in the variables
$\by = \set{\ya_{a,i}}{\alpha=1,\ldots,r;a\geq 1;i=1,\ldots,\ma_a}$.
Define the evaluation map
$\varphi_\mua$,
which sets all the variables in the same row of the (Young diagram
associated to the) partition $\mua$ to the same value,
$\xa_{a,i,j} \mapsto \ya_{a,i}$. The effect of the evaluation map
on the variables corresponding to the partition $\mua$ is shown in
figure \ref{evaluationmap2}. We define the evaluation map
$\varphi_\bmu : \overline\cF \rightarrow \cH$ to be 
$\varphi_\bmu = \prod_{\alpha=1}^{r} \varphi_\mua$.

By \eqref{intc}, $\varphi_\bmu (g(\bx)) = 0$ (where $g(\bx)$ is as in
\eqref{rationalf} with $\overline{G}(\bx) \in \overline{\cF}$), if any of the
partitions $\mua$ has a part which is greater than $k$. Hence, in the
following, we will assume that none of the partitions has a part
greater than $k$, and refer to these (multi)-partitions as $k$-restricted.

Our strategy will be to study the image of 
$\overline{\cF}$ under the evaluation map.

\begin{defn}
Let $\cH_\bmu$ be the space of functions in the variables
$\by$, and let $\overline{\cH}_\bmu \subset \cH_\bmu$ be the
subspace spanned by functions of the form
\begin{equation}
\label{h1}
H (\by) = H_\bmu (\by) h(\by) \ ,
\end{equation}
where $h (\by)$ is an arbitrary polynomial in $\by$, symmetric
under the exchange $y_{a,i}^{(\al)}\leftrightarrow
y_{a,j}^{(\al)}$, and
\begin{equation}
\label{h2}
\begin{split}
&H_\bmu (\by) =\\&
\prod_{\substack{\alpha=1,\ldots,r\\(a,i)>(b,j)}}
(\ya_{a,i}-\ya_{b,j})^{2 A_{ab}}
\prod_{\substack{\alpha=1,\ldots,r-1\\(a,i);(b,j)}}
(\ya_{a,i}-\yaa_{b,j})^{-A_{ab}} 
\prod_{(a,i)} (\yb_{a,i})^{\max(0,a-l)} \ .
\end{split}
\end{equation}
Here, $A_{ab}=\min(a,b)$ and 
$(a,i)\in I_k\times I_{m_\alpha}$ (where $I_m = \{1,\ldots,m\}$).
The ordering $(a,i)>(b,j)$ is defined
as follows. The index $i$ increases downwards, and we say that
$(a,i)>(b,j)$ if $a>b$, or, if $a=b$, when $i<j$.
\end{defn}

Let us define a lexicographic ordering on multi-partitions. That is,
the usual lexicographic ordering is taken on partitions $\mu^{(\alpha)}$,
and $\bnu > \bmu$ if
$\nua=\mua$ for all $\alpha<\gamma$ and $\nu^{(\gamma)}> \mu^{(\gamma)}$.

Let $\ker \varphi_\bmu$ be the kernel of the evaluation map
$\varphi_\bmu$ acting on $\overline\cF$. We can now define
the subspaces
\begin{equation}
\Gamma_{\boldsymbol\mu} = \bigcap_{\boldsymbol{\nu}>\boldsymbol{\mu}}
\ker{\varphi_{\boldsymbol{\nu}}} \quad , \quad
\Gamma'_{\boldsymbol\mu} = \bigcap_{\boldsymbol{\nu}\geq\boldsymbol{\mu}}
\ker{\varphi_{\boldsymbol{\nu}}} \ .
\end{equation}
Thus, $\Gamma_{\bmu}$ is the space of rational functions which are
annihilated by every evaluation map with $\bnu>\bmu$. 

By definition,
$\Gamma_{\bnu} \subset \Gamma_{\bmu}$ if $\bnu<\bmu$, and
$\Gamma'_{\bmu} \subset \Gamma_{\bmu}$. In addition, 
$\Gamma'_{(1^{m^{(1)}},\ldots,1^{m^{(r)}})}=\{0\}$.  Therefore,
$\Gamma_\bmu$ defines a filtration on $\overline\cF$.
Define the associated graded space
\begin{equation}
\Gr \Gamma = \bigoplus_{\bmu} \Gr_{\bmu} \Gamma \ ,
\end{equation}
where $\Gr_{\bmu} \Gamma = \Gamma_{\bmu} / \Gamma'_{\bmu}$ and the
sum is over multi-partitions of $\bm$.
The main purpose of this section is to prove
\begin{thm}
\label{thmigvs}
The induced map 
\begin{equation}
\label{emapgvs}
\overline{\varphi}_{\bmu} : \Gr_{\bmu} \Gamma \rightarrow \overline{\cH}_\bmu
\end{equation}
is an isomorphism of graded vector spaces.
\end{thm}
This is very similar to the proof found in \cite{FKLMM2} for the case
which corresponds to $\widehat{\sl}_3$, and we use the same ideas here.

To prove the theorem, we need to show three things. First, the
evaluation map
\begin{equation}
\label{emapg}
{\varphi}_\bmu : \Gamma_\bmu \rightarrow \overline{\cH}_\bmu 
\end{equation}
is well-defined. Second, it is surjective, and third, the induced map
\eqref{emapgvs} is well defined and injective.

\subsubsection{The evaluation map is well-defined}
To prove that the map
$\varphi_\bmu : \Gamma_\bmu \rightarrow
\overline{\cH}_\bmu$, is well 
defined, we must show that the rational
functions obtained after the evaluation are indeed of the form
\eqref{h1} and \eqref{h2}.
We will do this by showing that the structure of the poles and zeros
of the image of the functions \eqref{ratfb} in $\overline{\cF}_\bm$ under
the evaluation map is precisely of the form \eqref{h2}.

\begin{lemma} \label{lemintzeros}
Let $\overline{G}(\bx)\in \Gamma_\bmu$. Then, the function
$\varphi_\bmu (\overline{G}(\bx))$ has a zero of order at
least $2\min(a,a')$ when
$\ya_{a,i}=\ya_{a',i'} \ , \ \forall\alpha$.
\end{lemma}

\begin{proof}
The proof is independent of $\alpha$, and so we can use the argument used
in the case of $\widehat{\sl}_2$ in \cite{AKS}.
We will repeat that argument here for completeness.  

It is sufficient to consider 
the dependence of $\overline{G}(\bx)$ on the two sets of variables
of the same color $\alpha$, which we denote by
$\set{x_{a,i}}{i=1,\ldots,a}$ and $\set{x_{a',i'}}{i'=1,\ldots,a'}$.
We can assume that $a\geq a'$ without loss of generality.

We can carry out the evaluation map in two steps:
$\varphi_{\bmu}=\varphi^2\circ\varphi^1$. Here $\varphi^1$
consists of evaluating all the variables except the set
$\set{x_{a',i'}}{i'=1,\ldots,a'}$ and $\varphi^2$
consists of setting
$x_{a',1}=\cdots=x_{a',a'}=y_{a'}$ (note that under $\varphi^1$,
the variables $x_{a,1},\ldots,x_{a,a}$ are all set to $y_a$).

Let
 \begin{equation}
 g_1(y_{a};
 x_{a',1},\ldots,x_{a',a'})
 =\varphi^1 (\overline{G}(\bx)) \ .
 \end{equation}
Because $\overline{G} (\bx)\in \Gamma_\bmu$, 
$\overline{G}(\bx)$ is
annihilated by all $\varphi_{\bnu}$ with $\bnu>\bmu$. Therefore 
\begin{equation}
\left.
g_1(y_a;
 x_{a',1},\ldots,x_{a',a'})
\right|_{x_{a',i'}=y_{a}} =0 \ \hbox{for all $i'$},
\label{EQ:f_1_zero}
\end{equation}
because this corresponds to an evaluation corresponding to a multi-partition
greater than $\bmu$. Therefore,
\begin{equation}
g_1(y_{a}; x_{a',1},\ldots,x_{a',a'})=
\prod_{i'=1}^{a'} (x_{a}-x_{a',i'})
\tilde{g}_1(y_{a}; x_{a',1},\ldots,x_{a',a'}) \ .
\label{EQ:firstorder_zero}
\end{equation}
 
Now $g_1(y_{a}; x_{a',1},\ldots,x_{a',a'})$
was obtained from a symmetric function in $x_i^{(\al)}$, and so, for each $i'$,
\begin{equation}
\left.
\frac{\partial g_1}{\partial y_{a}}  
 \right|_{x_{a',i'}=y_{a}} =
a \left.\frac{\partial g_1}{\partial x_{a',i'}}
 \right|_{x_{a',i'}=y_{a}} \ .
 \label{EQ:derivative1}
\end{equation}
However (\ref{EQ:firstorder_zero}) tells us that, again for
each $i'$, 
\begin{equation} 
\left.
\frac{\partial g_1}{\partial y_{a}}
 \right|_{x_{a',i'}=y_{a}}  =
- \left.\frac{\partial g_1}{\partial x_{a',i'}}
 \right|_{x_{a',i'}=y_a} = 
\left.
\sideset{}{'}\prod_{i''=1}^{a'}
 (y_{a}-x_{a',i''})
 \tilde{g}_1\right|_{x_{a',i'}=y_{a}} \ ,
\label{EQ:derivative2}
\end{equation}
the prime on the product meaning that the term with $i''=i'$ is
to be omitted.
The only way to reconcile (\ref{EQ:derivative1}) with
(\ref{EQ:derivative2}) is for 
$\tilde g_1|_{x_{a',i'}=y_{a}}$ to be zero.
Thus the  zero at $x_{a',i'}=y_{a}$
is at least of order two    
\begin{equation}
g_1(y_{a}; x_{a',1},\ldots,x_{a',a'})=
 \prod_{i'=1}^{a'}
 (y_{a}-x_{a',i'})^2
\tilde{g}_2
(u_{a}; x_{a',1},\ldots,x_{a',a'}) \ .
\label{EQ:f_1_zero2}
\end{equation}
We now evaluate the right-hand-side of (\ref{EQ:f_1_zero2})
at $x_{a',1}=\cdots=x_{a',a'}=y_{a'}$
and, recalling the condition that $a\geq a'$, we have
\begin{equation}
\varphi_\bmu(G(\bx)) =
\prod (y_{a}-y_{a'})^{2 A_{a,a'}} \widetilde{G} \ .
\label{EQ:A_matrix}
\end{equation}
\end{proof}

\begin{lemma} \label{lemserpoles}
The image under the evaluation map $\varphi_\bmu$
of any function in $\overline{\cF}$ (and hence $\Gamma_\bmu$)
has a pole of maximal order $\min(a,a')$ whenever
$\ya_{a,i}=\yaa_{a',i'}$.
\end{lemma}

\begin{proof}
We will prove this lemma by looking
at the zeros of ${g}(\bx)$, which arise because we need to
satisfy the Serre relations,
${g}|_{\xa_{1}=\xa_{2}=\xaa_{1}}=0$ and 
${g}|_{\xa_{1}=\xaa_{1}=\xaa_{2}}=0$ for $\alpha=1,\ldots,r-1$.
These relations depend on two sets of variables only. 

Consider the dependence of ${g}$ on the two sets of variables
$x_i = \xa_i$, with $i=1,\ldots,a$ and $\bar{x}_j = x^{(\al\pm1)}_j$,
with $j=1,\ldots,a'$. Under the evaluation map,
these variables map to $\varphi_\bmu (x_i) = y$ and 
$\varphi_\bmu (\bar{x}_i) = \bar{y}$ respectively. 

Note that $x$ and $\bar{x}$ are variables corresponding to two
adjacent roots.  Again without loss of generality, assume that $a \geq a'$.

When $x_1=\bar{x}_1=\bar{x}_j$ or $x_1=x_j=\bar{x}_1$,
${g}$ vanishes, so we find
\begin{equation}
\begin{split}
&{g} (x_1,\ldots,x_a;\bar{x}_1,\ldots,\bar{x}_{a'};\ldots)
|_{x_1=\bar{x}_1=z_1} = \\
&\prod_{i=2}^{a} (x_i-z_1) \prod_{j=2}^{a'} (\bar{x}_i-z_1)
{g}' (z_1;x_2,\ldots,x_a;\bar{x}_2,\ldots,\bar{x}_{a'};\ldots) 
\end{split}
\end{equation}
Repeating the argument for ${g}'$ we find
\begin{equation}
\begin{split}
&{g}' (x_2,\ldots,x_a;\bar{x}_2,\ldots,\bar{x}_{a'};\ldots)
|_{x_2=\bar{x}_2=z_2} = \\
&\prod_{i=3}^{a} (x_i-z_2) \prod_{j=3}^{a'} (\bar{x}_i-z_2)
{g}'' (z_1,z_2;x_3,\ldots,x_a;\bar{x}_3,\ldots,\bar{x}_{a'};\ldots) 
\end{split}
\end{equation}
We can repeat this argument $a'$ times with the result
\begin{equation} \label{serarg}
\begin{split}
&{g} (x_1,\ldots,x_a;\bar{x}_1,\ldots,\bar{x}_{a'};\ldots)
|_{\{x_i=\bar{x}_i=z_i\}_{i=1}^{a'}} = \\
&\prod_{i=1}^{a'} \prod_{j=i+1}^{a} (x_j-z_i)
\prod_{i=1}^{a'} \prod_{j'=i+1}^{a'} (\bar{x}_{j'}-z_i)
\tilde{g} (z_1,\ldots,z_{a'};x_{a'+1},\ldots,x_a;\ldots) \ .
\end{split}
\end{equation}
We find that $\varphi_\bmu ({g})$ has a zero of order at least
$a a' -\min(a,a')$ when $y=\bar{y}$, by counting the number of zeros in
\eqref{serarg} and using that $a'\leq a$.
Taking into account the poles of \eqref{ratfb}, which
after applying the evaluation map becomes a pole of order $a a'$ when
$y=\bar{y}$, 
we find that the image of $\overline{\cF}_\bm$ has a pole or order at
most $\min(a,a')$, when $\xa_{a,j}=x^{(\al\pm 1)}_{a',j'}$.
\end{proof}

\begin{lemma} \label{lemirrpoles}
The image of $\varphi_{\bmu}$ acting on a function
$\overline{G} \in \Gamma_\bmu$ has a zero of order at least
$\max(0,a-l)$ when $\yb_{a,i}=0$.
\end{lemma}
\begin{proof}
To prove this lemma, we will study the 
effect of the evaluation map on ${g}(\bx)$ in eq. \eqref{rationalf}. We
focus on the variables of a row of length $a$ (where we assume that
$a>l$), $\set{\xb_j}{j=1,\ldots,a}$. Under the evaluation map, these
variables map to $\varphi_\bmu ( \xb_j ) = \yb$.

We know that the function
\begin{equation}
g_1(\xb_1,\ldots,\xb_a) =
g (\bx)|_{\xb_1=\cdots=\xb_{l}=0}
\end{equation}
contains a factor $\prod_{j=l+1}^{a} \xb_j$, because it
vanishes if any of the remaining variables $\xb_j$ is set to
zero (because of the condition \eqref{irrepc} on ${g}(\bx)$).
Thus, the image of $g_1$
under the evaluation map has a zero of order at least
$\max (0,a-l)$ whenever $\yb_{a,i}=0$.
\end{proof}

\begin{lemma} \label{wdef}
The map $\varphi_\bmu : \Gamma_\bmu \rightarrow \overline{\cH}_\bmu$
is well defined.
\end{lemma}
\begin{proof}
This follows from the lemmas
\ref{lemintzeros}, \ref{lemserpoles}, \ref{lemirrpoles}
and the definition of the space $\overline{\cH}_\bmu$.
\end{proof}

\subsubsection{Proof of surjectivity}
We will continue with the proof that
the map \eqref{emapg} is surjective. We have to prove that
for each function of the form defined by \eqref{h1} and \eqref{h2},
there is at least one function in the 
pre-image in $\Gamma_\bmu$. We do this by
explicitly giving the form of these pre-images, showing that they
are elements of $\overline{\cF}$ and finally, proving that
these pre-images are indeed in the kernel of $\varphi_\bnu$ for
each $\bnu>\bmu$, which shows that they are in $\Gamma_\bmu$.

For each ($k$-restricted) multi-partition $\bmu$, we consider the
function
\begin{equation} \label{gpreim}
F(\bx) = \frac{\sym f(\bx)}{p(\bx)} \ ,
\end{equation}  
where $f(\bx)$ and $p(\bx)$ are a polynomials of the form
(we identify the variables $\xa_{a,i,a+1}=\xa_{a,i,1}$)
\begin{align}
\label{g2preimage}
f (\bx) &=  \tilde{f} (\bx) 
\prod_{\substack{\alpha,a,i\\j>l}} \xb_{a,i,j}
\prod_{\substack{\alpha\\a,i,j\\a',i';j'\neq j}}
(\xa_{a,i,j}-\xaa_{a',i',j'}) 
\prod_{\substack{\alpha\\(a,i)>(a',i')\\j=1,\ldots,\ma_{a'}}}
(\xa_{a,i,j}-\xa_{a',i',j})
(\xa_{a,i,j+1}-\xa_{a',i',j}) \\
\label{preimpoles}
p(\bx) &= 
\prod_{\substack{\alpha=1,\ldots,r-1\\a,i,j\\a',i',j'}}
(\xa_{a,i,j}-\xaa_{a',i',j'}) \ ,
\end{align}
where $\tilde{f} (\bx)$ is an arbitrary polynomial. 
The symmetrization is over each of the $r$ sets of variables
$\{ \xa_i \}$ with the same value of $\alpha$.
As we did before, we will drop as many indices as possible in the
following lemmas.

\begin{lemma} \label{ginf}
The functions $F(\bx)$ of \eqref{gpreim} are elements of
$\overline{\cF}$
\end{lemma}
\begin{proof}
We have to show that $f(\bx)$
satisfies the vanishing conditions \eqref{intc}, \eqref{serrec} and
\eqref{irrepc}. 
First of all, we easily see that $f(\bx)$ is zero when any $k+1$
variables of the same color are set
to the same value. Because the partitions have rows of maximum length
$k$, these $k+1$ variables can not all be placed in the same row, which
implies that the factor
$\prod (\xa_{a,i,j}-\xa_{a',i',j'})$ evaluates to zero under
$\varphi_\bmu$.

To show that the
Serre relations are satisfied, we have to show that the zeros
\begin{equation} \label{serzeros}
\prod_{\substack{\alpha\\a,i,j\\a',i';j'\neq j}}
(\xa_{a,i,j}-\xaa_{a',i',j'})
\end{equation}
satisfy the Serre relations. Let $x_{a,j}= \xa_{a,i,j}$ and
$\bar{x}_{a',j'} = \xaa_{a',i',j'}$,
for some choice of $\alpha$, $i$ and $i'$.

For every $x$, there is a zero with every $\bar{x}$,
except those appearing in the column which has the same number as the
$x$ (i.e. for $j=j'$). Note that if we set two variables
$x$, which belong to the same column, to the same value, 
$f(\bx)$ is zero, because the factor
$\prod (\xa_{a,i,j}-\xa_{a',i',j'})$ is zero in that case.
Hence, we set
$x_{a,j}=\bar{x}_{a',j'}=\tilde{x} \ (j'\neq j)$. Focusing
on this variable, we find the following zeros
$(\tilde{x}-\bar{x}_{a,i})(\tilde{x}-\bar{x}_{a',i'})
\prod_{a'';i''\neq i,i'}(\tilde{x}-\bar{x}_{a'',i''})^2$
So, indeed $\tilde{x}$ has zero with every $\bar{x}$. Similarly, we find
that there is at least a zero of order one when we set
$x_1=\bar{x}_1=\bar{x}_2$.

To complete the proof of this lemma, we need to show that
$f(\bx)$ satisfies the condition \eqref{irrepc}. This easily follows
form the factor $\prod_{j> l} \xb_{a,i,j}$, combined with the
zeros which give rise to the condition \eqref{intc}.
\end{proof}

\begin{figure}[ht]
\begin{center}
\includegraphics{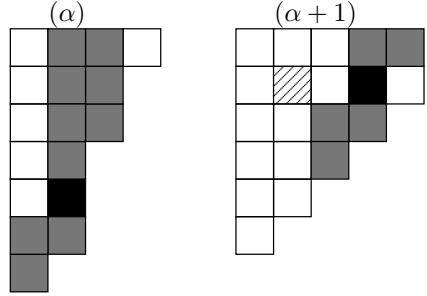}
\end{center}
\caption{A violation of the Serre relations if the zero corresponding
to the black squares is removed from \protect\eqref{g2preimage}.
The left partition corresponds to the
variables of color $(\alpha)$, the right one to color $(\alpha+1)$.
The `slanted' box is the third
variable, in addition to the two black ones, for which the Serre
condition is violated. The gray boxes denote the zeros with the
variable corresponding to the black box of the same partition, coming
from the integrability conditions.}
\label{serrev2}
\end{figure}

\begin{remark}
It is instructive to note that all the zeros in \eqref{serzeros}
are necessary to satisfy the Serre relations. We need to show that
if we remove any of these zeros, we will violate a Serre
relation.
\end{remark}
To show that this is true,
it is important that we take the zeros between variables
of the same color into account. Let us remove the zero
$(x_{a,j}-\bar{x}_{a',j'})$,  where $j \neq j'$. Without loss of generality,
we can assume that $j<j'$.
The two variables are indicated in figure \ref{serrev2} by the black boxes.
The gray boxes denote the zeros with the variables corresponding to
the black box from the same partition.
All we need to do is show that there is at least one
variable, of either partition, such that when this variable is set to
the same value as the two `black variables', we do not get a zero,
and thus violate a Serre relation. This variable is taken to be of
color $(\alpha+1)$, (if $j > j'$, it is of color $(\alpha)$).
More precisely, it is the variable $\bar{x}_{a',j}$,
taken from the same row as $\bar{x}_{a',j'}$ (denoted by the `slanted' box),
which always exists, because $j<j'$.

There is no zero at
$\bar{x}_{a',j'}=\bar{x}_{a',j}$, because both
variables are taken from the same row. In addition, there is no zero
at $x_{a,j}=\bar{x}_{a',j}$, because it is not
present in the factor \eqref{serzeros} and the zero at
$x_{a,j}=\bar{x}_{a',j'}$
is the one we removed. We conclude that after we remove the
(arbitrary) zero at $x_{a,j}=\bar{x}_{a',j'}$, we do not have a zero when
$x_{a,j}=\bar{x}_{a',j'}=\bar{x}_{a',j}$. Thus, we have shown that by
removing any of the zeros in \eqref{serzeros}, we violate a Serre
condition. We conclude that the zeros are indeed necessary.

\begin{lemma} \label{ginkernel}
The function $F(\bx)$ of \eqref{gpreim} associated to a
$k$-restricted multi-partition
$\bmu$ is an element of the kernel of $\varphi_\bnu$ for
any $\bnu>\bmu$.
\end{lemma}
\begin{proof}
Let us take a $\bnu>\bmu$, and let $\nua$ be the
first partition such that $\nua>\mua$. We will focus on
the variables $\xa$ and show that the function
$F(\bx)$ can not be non-zero under the evaluation map $\varphi_\bnu$.

Two variables in the same column of $\mua$ have a zero, so they
can not be placed in the same row in $\nua$, if the result is to be
non-zero, because in that case, acting with the evaluation map
gives a zero.

However, because
$\nua>\mua$, we can not avoid placing variables of the same
column in $\mua$ in the same row of $\nua$. To show this,
let us denote the length of the rows of the partitions by
$\mu^{(\alpha)}_i$ and $\nu^{(\alpha)}_i$, such that the index
$i$ is increasing
going downwards. The only way to avoid placing variables of the same
column of $\mua$ in the same row of $\nua$ is by placing
the variables of $\mua$ in a rows of the same length in
$\nua$. However, because $\nua>\mua$, there will be an
$\tilde{\imath}$ such that
$\nu^{(\alpha)}_{\tilde{\imath}}>\mu^{(\alpha)}_{\tilde{\imath}}$.
Let us focus on the smallest $\tilde{\imath}$. We have to place a
variable of a row $\mu^{(\al)}_i$ with $i>\tilde{\imath}$ in the row
$\nu^{(\al)}_{\tilde{i}}$. Because
$\mu^{(\alpha)}_i \leq \mu^{(\alpha)}_{\tilde{\imath}}$,
this variable belongs to
the same column of another variable in $\nu^{(\alpha)}_{\tilde{\imath}}$.
We conclude that $F(\bx)$ is zero under the evaluation map
$\varphi_\bnu$ with $\bnu>\bmu$.
\end{proof}

\begin{lemma} \label{gingamma}
The function $F(\bx)$ of \eqref{gpreim} is an element of
$\Gamma_\bmu$.
\end{lemma}
\begin{proof}
This follows from the lemmas \ref{ginf} and \ref{ginkernel}
\end{proof}

As a last step in the proof of surjectivity, we have to show that
the image of $F(\bx)$ under the evaluation map is indeed of the form
\eqref{h1} and \eqref{h2}.
In particular, it contains as a factor the functions
$h(\by)$, which are symmetric under the exchange of variables
$\ya_{a,i} \leftrightarrow \ya_{a,i'}$.  

\begin{lemma} \label{evalgsmh}
The image of $F(\bx)$ under the evaluation map
$\varphi_{\bmu}$ is a scalar multiple of the function $H(\by)$ in
\eqref{h1}.
\end{lemma}
\begin{proof}
To prove this lemma, we can follow the same approach as we did
in our paper on the $\widehat{\mathfrak{sl}}_2$ case, because
the argument does not depend on the color of the variables.
We will focus on the variables $\xa$, and determine the permutations
$\sigma$, for which $\varphi_\bmu(f(\sigma(\xa)))$ is non-zero.
So, we consider
\begin{equation}
\sum_{\sigma\in \cS_{m^{(\alpha)}}} f( \sigma\{\xa\}) \ .
\end{equation}
In the following, we will omit the label $\alpha$.
Recall that the variable $x_{a,i,j}$ corresponds to the
$j^{th}$ column in the $i^{th}$ row of length $a$. Under the
evaluation map, $x_{a,i,j}\mapsto y_{a,i} \; \forall j$.

Suppose that for some $\sigma$, we have
$\sigma (x_{a,i,j})=x_{a',i',j'}$ with $(a',i')<(a,i)$
and that $(a,i)$ is the largest row for which this is true. This
means that all rows above $(a,i)$ undergo only a permutation within
the row. Suppose that the pre-factor 
\begin{equation}
\varphi_\bmu \circ \sigma \left(
\prod_{(a,i)>(a',i')}
(x_{a,i,j}-x_{a',i',j})
(x_{a,i,j+1}-x_{a',i',j})
\right) 
\end{equation}
is to be non-zero. Then $x_{a',i',j'}$ can not be in  a
column directly below or to the left of the
permutation image of any other element from row $(a,i)$.
This means that at least one other element from row
$(a,i)$ should be mapped to a row below $(a,i)$. If
it is mapped to the row $(a',i')$ it can appear in any
column other than $j'$. If it mapped to any other row, it can
appear in any other column  than $j'$ and an adjacent column
(to the right or left depending on whether it is above or
below $(a',i')$.) Now we repeat this argument for this
new element, concluding that at least one more element of
row $(a,i)$ is mapped to a lower row, and so forth,
until eventually we find that all elements are permuted
to a row below $(a,i)$. 
If the elements are permuted to the same row, they can be placed
adjacent columns. Elements which are permuted to different rows can
not be placed in adjacent columns, this being due to the factor
linking adjacent columns in the pre-factor.  
There are at most $a$ columns in $\mu^{(\alpha)}$ in rows below
$(a,i)$,
and hence the elements must all appear in the same row, which is
therefore of length $a$. Thus all the variables in rows of length $a$
are mapped to another row of length $a$, for the same
reason. As a result, the only permutations which give a non-zero
contribution to $\varphi_\bmu (f(\sigma(\xa)))$ are those that
permute variables within each row, or those that
permute rows of equal length. Under the evaluation map, the
former contribute equal terms to the sum, while row
interchanges correspond to the symmetrization over
the variables $\ya_{a,i}$ with the same values of $\alpha$ and
$a$ in $h(\by)$. Note that the other factors in the
function $F$ are symmetric under the permutation of rows of equal
length, so these factors do not interfere with the argument above.
\end{proof}

\begin{lemma} \label{surjective}
The map $\varphi_{\bmu} : \Gamma_\bmu \rightarrow \overline{\cH}_\bmu$
is surjective.
\end{lemma}
\begin{proof}
This follows from the lemmas \ref{gingamma} and \ref{evalgsmh}.
\end{proof}

\subsubsection{Injectivity proof}
\begin{lemma} \label{injective}
The induced map
$\overline{\varphi}_{\bmu} : \Gr_{\bmu} \Gamma \rightarrow \overline{\cH}_\bmu$
\eqref{emapgvs} is well defined and injective.
\end{lemma}
\begin{proof}
To prove that the map \eqref{emapgvs} is well defined, we use 
lemma \ref{wdef} and observe that
the image of $\Gamma'_{\bmu}$ under $\varphi_{\bmu}$ is zero by using
the definition of $\Gamma'_{\bmu}$. It follows that we can
define the induced map $\overline{\varphi}_\bmu$  acting on the quotient
$\Gr_{\bmu} \Gamma= \Gamma_{\bmu}/\Gamma'_{\bmu}$. Moreover, the difference
between two different functions in $\Gamma_{\bmu}$ that map to the same
rational function in $\overline{\cH}_\bmu$ is in $\Gamma_{\bmu}$. Hence, the
map is also injective. 
\end{proof}

We have now completed the proof of theorem \ref{thmigvs},
because the theorem follows from the lemmas \ref{injective} and
\ref{surjective}.

\bigskip
The map \eqref{emapgvs} is degree preserving, and thus we can count
the functions of homogeneous degree $d$ 
in $\overline{\cH}_\bmu$ to obtain the character of the space
$\overline{\cF}$.

To compute the character of $\cF_\lambda$, we add the poles $\prod
(\xb_{a,i,j})^{-1}$, which are present in the functions $G(\bx)$ in
\eqref{rationalf}. The only thing in the calculation of the character
which changes is the fact that due to these poles, the zeros $\prod
(\yb_{a,i})^{\max(0,a-l)}$ in \eqref{h2} become poles $\prod
(\yb_{a,i})^{-\min(a,l)}$.
 
\subsection{Character of the dual space}
\label{cds}

Using the results of the previous section, we can calculate the
character of the dual space $\cF_\lambda$, where
$\lambda = l \omega_\beta$.

First, let us define the character of $W_\lambda$ as follows:
\begin{equation}\label{wchardef}
\ch_{q} W_\lambda =  \sum_{d,m^{(\alpha)}}   \dim
W_\lambda[\mathbf m]_d \
q^d e^{\lambda - \bomega^T C_r \bm },
\end{equation}
where $W_\lambda[\mathbf m]_d$ is the subspace generated by elements
in $U(\wn)$ of homogeneous degree $m^{(\alpha)}$ in $f_\alpha$,
and homogeneous degree $-d$ in $t$. Here,
$\bomega=(\omega_1,\ldots,\omega_r)^T$.

The space $\mathcal F_\lambda$ is a space of functions in the
variables $x_i^{(\alpha)}$. If we define its $(\mathbf m,d)$-graded
component to be the space of functions in $m^{(\alpha)}$ variables
$x_i^{(\alpha)}$ and total homogeneous degree $\widetilde{d}$ in all
the variables, then, due to the way we defined the generating
functions $f_\alpha(x)$ (or, equivalently, the coupling), we have that
$\mathcal F_\lambda[\mathbf m]_{\widetilde{d}}$ is the dual to
$W_\lambda[\mathbf m]_d$ where $d = \widetilde{d} + \sum_\alpha
m^{(\alpha)}$.

Thus,
\begin{equation}
\ch_q W_\lambda = \sum_{\bm} \ch W_{\lambda}[\bm]
=\sum_\bm \sum_d q^{d+\sum_{\alpha} m^{(\alpha)}}
e^{\lambda-\bomega^T C_r \bm}
\dim
(\cF_{\lambda}[\bm])_d \ ,
\end{equation}
where $\dim (\cF_{\lambda}[\bm])_d$ denotes the dimension of the
subspace of functions in $\cF_{\lambda}[\bm]$ which have homogeneous
degree $d$. The powers of $z$ correspond to the components of the
weights in terms of the simple roots. Recall that here, $\lambda =
\lambda_\beta = l\omega_\beta$.

We will calculate this character by actually summing over all the
functions in $\cH$, and counting their homogeneous degree.  The
character of the space of symmetric functions $h (\by)$ in $m_a^{(\al)}$
variables is given by
\begin{equation}
\frac{1}{\prod_{\alpha=1}^{r} \prod_{a=1}^{k} (q)_{\ma_a}} \ ,
\end{equation}
where $(q)_m = \prod_{i=1}^m (1-q^i)$ for $m\in\NN$ and $(q)_0=1$.

The homogeneous degree of the rational function $H_\bmu (\by)$,
combined with the additional poles ${\prod_{(a,i)} (\ya_{a,i})^{-a}}$
is given by
\begin{equation}
\deg \left(\frac{H_\bmu (\by)}{\prod_{(a,i)} (\ya_{a,i})^{a}}
\right) = \sum_{\alpha,\alpha',a,a'}
\frac{1}{2} m_{a}^{(\alpha)} (C_r)_{\alpha,\alpha'} A_{a,a'}
m_{a'}^{(\alpha')}
- \sum_{a} A_{a,l} m_{a}^{(\beta)}
- \sum_{\alpha} m^{(\alpha)} \ .
\end{equation}
It follows
that the character of $W_\lambda^{(0)}$ is
\begin{equation}
\label{charw0}
\ch_q W_{l\omega_\beta}^{(0)} =
\sum_{\vbm \in \Z_{\geq 0}^{r\times k}}
\frac{
q^{\frac{1}{2} {\vbm}^T (C_r\otimes A){\vbm}
- (\id\otimes A \vbm)_{l}^{(\beta)}}}
{(q)_{\vbm}}e^{l\omega_\beta-\bomega^T C_r\bm}
.
\end{equation}
Here $(A)_{a,b}=\min(a,b)$ is a $k\times k$ matrix, and $C_r$ is the
Cartan matrix of $\sl_{r+1}$. Also, $\vbm$ denotes the vector
$(m_1^{(1)},\ldots,m_k^{(1)};\cdots;m_1^{(r)},\ldots,m_k^{(r)})$.
We made use of the definition
$$
(q)_{\vbm}={\prod_{\alpha=1}^{r}\prod_{a=1}^{k} (q)_{m_{a}^{(\alpha)}}} \
 .
$$

\section{Characters for rectangular highest-weight
$\widehat{\sl}_{r+1}$-modules}
\label{awt}
In this section, we will show that we can use the characters of the
principal subspace $W_\lambda$ to obtain the character of the full
integrable module $V_\lambda$. We will be able to do this by using the
invariance of the weight multiplicities of $V_\lambda$ under the
action of the affine Weyl group, in particular the affine Weyl
translations $t_\alpha$. More specifically, we will show that acting
with an affine Weyl translation on the principal subspace, and taking
an appropriate limit, we obtain the full integrable module.

Let $\Lambda$ be an affine weight of level $k$. It can be written as
$$\Lambda = \lambda + k \Lambda_0 - m \delta,$$ where $\lambda$ is the
weight with respect to $\h\in\sl_{r+1}$. Let $t_{\alpha}$ be
the affine Weyl translation corresponding to the
root $\alpha$ (see \cite{Kac}, equation (6.5.2)),
 and define the translation $t_\bN = \prod_i t_{N_i\alpha_i}$,
where $\bN = (N_1,\ldots,N_r)^T$. Then
\begin{equation}
\label{awtsr}
t_{\bN}(\Lambda) =
\lambda + k \bN^T\cdot\boldsymbol\alpha + k\Lambda_0 -
(m+\bN^T\cdot\mathbf l+ \frac{1}{2} k\bN^T C_r\bN)\delta.
\end{equation}
Again, $\mathbf l = (l_1,\ldots,l_r)^T$ where $\lambda = \sum_i
l_i\omega_i$, and $\boldsymbol\alpha = (\alpha_1,\ldots,\alpha_r)^T$.
Also note that $\bfalpha$ in terms of the weights is given by
$\bfalpha = C_r \bomega$.

Consider the principal subspace $W^{(\bN)} = U(\wn) t_{\bN}
v_\lambda$. It has a dual space description which is similar to
$\cF_\lambda$, if we choose the vector $\bN$ carefully.
Given that if $f_\al[m] v_\lambda = 0$, then $f_\al[m+(C_r\cdot\bN)_\al]
t_{\bN}v_\lambda = 0$ (since the Weyl group preserves weight space
multiplicities), we choose $\bN$ such that $(C_r\cdot\bN)_\al =
2N$ for all $\al$, for some $N\in\Z_+$. In the case of $\sl_{r+1}$, we
have $(\bN)_i= N i (r+1-i)$. 

Then $f_\al[2N+\delta_{\al,\beta}] t_\bN 
v_\lambda = 0$,  where $\lambda = l\omega_\beta$, and
$f_\al[2N-1+\delta_{\al,\beta}]t_\bN v_\lambda\neq 0$.

Note also that the extremal vector $t_{\bN} v_\lambda$ is a basis for the
one-dimensional weight subspace of weight 
$$t_{\bN}(\lambda ) = \lambda + k \bN^T\boldsymbol\alpha + k \Lambda_0
- \bigl((\lambda,\bN^T\boldsymbol\alpha)+\frac{1}{2}k\bN^T
C_r\bN\bigr)\delta.
$$
In the case of interest here, this becomes
$$
t_{\bN} (l\omega_\beta) = l\omega_\beta + k \bN^T\boldsymbol\al +
k\Lambda_0 -(l N_\beta + k N |\bN|)\delta \ ,
$$
where $|\bN| = \sum_i N_i$.

Thus, the space dual to $W^{(\bN)}_\lambda$ is the space of functions
of the form
$$
\prod_{\alpha,i}(x_i^{(\al)})^{-2N} G(\bx),
$$
where $G(\bx)$ is the function in equation (\ref{rationalf}). 

Thus, we find that the character of $W^{(\bN)}_{l \omega_\beta}$
differs from the character of $W^{(0)}_{l \omega_\beta}$ by a change
in the exponent of $q$ by 
$l N_\beta + k N |\bN| -2N|\bm|$ (where $|\bm| = \sum_\alpha \ma$)
and a change in the weight by $\bomega^T C_r k \bN$, which leads to 
$$
\ch_{q} W_{l\omega_\beta}^{(\bN)} = \sum_{\vbm\in \Z_{\geq 0}^{r\times k}}
\frac{q^{\half \vbm^T (C_r\otimes A)\vbm - (\id\otimes A\vbm)_l^{(\beta)}}}
{(q)_\vbm}
q^{l N_\beta + k N |\bN| -2N|\bm|}
e^{l\omega_\beta -\bomega^T C_r(\bm-k\bN)}.
$$

A form suitable for taking the limit $N\to\infty$ is obtained by
eliminating the summation variable $m_k^{(\alpha)}$ in favor of
$m^{(\alpha)} = \sum_{a=1}^k a m_a^{(\al)}$. 
This gives for the character of $W^{(0)}_{l\omega_\beta}$
(we define $\overline{m}^{(\alpha)} = \sum_{a=1}^{k-1} a \ma_a$)
\begin{equation}
\label{charw0trans}
\begin{split}
&\ch_{q} W^{(0)}_{l \omega_\beta} =
\sum_{\bm \in \Z_{\geq 0}^{r} }
q^{\frac{1}{2k} \bm^T C_r \bm
- \frac{1}{k}l m^{(\beta)}}
e^{l\omega_\beta-\bomega^T C_r\bm}
\\&\times
\sideset{}{'}\sum_{\vbm \in \Z_{\geq 0}^{r\times(k-1)}}
\frac{q^{\frac{1}{2} \vbm^T (C_r\otimes C^{-1}_{k-1})\vbm
- \bigl((\id\otimes C^{-1}_{k-1})\vbm\bigr)_l^{(\beta)} \delta_{l<k}}}
{\prod_{\alpha=1}^{r}\prod_{a=1}^{k-1}
(q)_{\ma_a} (q)_{\frac{\ma-\overline{m}^{(\alpha)}}{k}}} \ ,
\end{split}
\end{equation}
where the prime on the sum denotes the constraints
$\overline{m}^{(\al)} \leq \ma$ and $\overline{m}^{(\al)}\equiv \ma
\mod k$. Here,
$C_{k-1}$ denotes the Cartan matrix of $\sl_k$.
The symbol $\delta_{l<k}$ is $1$ for the integers in the range
$l=1,\ldots,k-1$ and zero otherwise.

The character of $W_{l\omega_\beta}^{(\bN)}$ has an extra factor of 
$$
q^{l N_\beta + \half k\bN^T C_r\bN -\bm^T C_r\bN}
$$
Combining this power of $q$ with the power in the first line of 
(\ref{charw0trans}), we use the change of variables
$\widetilde{m}^{(\alpha)} = m^{(\alpha)}-kN_\alpha$, since the
combined power in this new variable is
\begin{eqnarray*}
{\frac{1}{2k} \bm^T C_r \bm + \half k \bN^T C_r \bN - \bm^T C_r\bN
- \frac{l}{k}(m^{(\beta)}-k N^{(\beta)})}
 = \frac{1}{2k} \widetilde{\bm}^T C_r \widetilde{\bm} - 
\frac{l}{k}\widetilde{m}^{(\beta)}.
\end{eqnarray*}

Making this substitution, we have
\begin{equation}
\begin{split}
\label{wnchar}
&\ch_q W^{(\bN)}_{l\omega_\beta} =\\&
\sum_{\widetilde{\bm} \geq - k\bN}
{q^{\frac{1}{2k} \widetilde{\bm}^T C_r\widetilde{\bm}
- \frac{1}{k} l\widetilde{m}_\beta}}
e^{l\omega_\beta -\bomega^T C_r\widetilde{\bm}}
\sideset{}{'}\sum_{\vbm\in\Z_{\geq 0}^{r\times(k-1)}}
\frac{q^{\frac{1}{2} \vbm^T (C_r\otimes C^{-1}_{k-1})
\vbm - \bigl((\id\otimes C^{-1}_{k-1})\vbm\bigr)^{(\beta)}_l \delta_{l<k}}}
{\prod_{\alpha=1}^{r}\prod_{a=1}^{k-1}
(q)_{\ma_a} (q)_{\frac{\widetilde{m}^{(\al)}-\overline{m}^{(\al)}}{k}+ N_\alpha}} \ ,
\end{split}
\end{equation}
where the prime denotes the constraints
$\overline{m}^{(\al)} \leq \widetilde{m}^{(\al)} + k N_\al$ and
$\overline{m}^{(\al)} = \widetilde{m}^{(\al)} \mmod k$. 

We can now easily obtain the characters of the integrable level-$k$
modules corresponding to rectangular highest weights
by taking the limit $N\to\infty$ while keeping $\widetilde{\bm}$
finite. This gives 
\begin{equation}
\label{vchar}
\begin{split}
&\ch_q W^{(\infty)}_{l\omega_\beta} = 
\sum_{\widetilde{\bm}\in\Z^{r}}
q^{\frac{1}{2k} \widetilde{\bm}^T C_r\widetilde{\bm}
- \frac{1}{k}l \widetilde{m}^{(\beta)}} e^{l\omega_\beta -\bomega^T
C_r\widetilde{\bm}} \\
&\qquad \times \frac{1}{(q)^{r}_\infty}
\sideset{}{'}\sum_{\vbm\in\Z_{\geq 0}^{r\times(k-1)}}
\frac{q^{\frac{1}{2} \vbm^T (C_r\otimes C^{-1}_{k-1})\vbm
- \bigl((\id\otimes C^{-1}_{k-1})\vbm\bigr)^{(\beta)}_l \delta_{l<k}}}
{(q)_{\vbm}}\ ,
\end{split}
\end{equation}
with the constraint
$\overline{m}^{(\al)} = \widetilde{m}^{(\al)} \mmod k$.

The nice feature of this character formula is that it manifestly
splits the character into a sum over all the finite weights, each of
which contributes a string function to the full character. These string
functions are proportional to `the second line' of equation
\ref{vchar}.

We can make the appearance of the characters slightly more
compact, by rewriting it in terms of the $r\times k$-vector $\vbm$ again.
This results in
\begin{equation}
\label{chwinfrect}
\ch_q W^{(\infty)}_{l\omega_\beta} = \frac{1}{(q)_\infty^{r}}
\sum_{\substack{\vbm\\m_k^{(\alpha)}\in\Z,m_{a<k}^{(\alpha)}\in\Z_{\geq 0}}}
\frac{
q^{\frac{1}{2} \vbm^T (C_r\otimes A)\vbm
- \bigl((\id\otimes A)\vbm\bigr)_l^{(\beta)}}}
{\prod_{\al=1}^{r}\prod_{a<k}(q)_{m_a^{(\al)}}}
e^{l\omega_\beta-\bomega^T C_r\bm} \ ,
\end{equation}
which, again, holds in the case of rectangular
representations. Comparing this to the known character formul\ae\ for
the integrable representations which appear, for example, in
\cite{Ge2}, we see that this is indeed the character of the
integrable, level-$k$ $\widehat{\sl}_{r+1}$-module with highest weight
$\lambda = l\omega_\beta$, i.e.\ the module $V_{l\omega_\beta}$.
Hence, we have the following result
\begin{thm}
The character of the integrable, level-$k$ $\widehat{\sl}_{r+1}$-module
with highest weight $\lambda=l\omega_\beta$ is given by
$$
\ch_q V_{l\omega_\beta} = \ch_q W^{(\infty)}_{l\omega_\beta} \ ,
$$
where $\ch_q W^{(\infty)}_{l\omega_\beta}$ is given by equation
\eqref{chwinfrect}.
\end{thm}

The remainder of the paper will be devoted to obtaining character
formul\ae\ for general irreducible representations.

\section{Conformal blocks and their dual spaces}
\label{confb}
\subsection{Modules localized at $\zeta\neq 0$}

Above, we considered the standard action of the central extension of
the loop algebra, $\widetilde{\g}=\g\otimes \C[t,t^{-1}]$ on
integrable modules $V_\lambda$ of level $k$. Such modules can be
considered as ``localized'' at the point $0$. 

For a generic point $\zeta\in\C P^1$, Let $t_\zeta=t-\zeta$ denote a
local variable at $\zeta$, and consider the action of the current
algebra $\widetilde{\g}_{(\zeta)} = \g \otimes
\C[t_\zeta,t_\zeta^{-1}]$ on a module $V_\lambda(\zeta)$, ``localized''
at the point $\zeta$, which is isomorphic to $V_\lambda$.

Specifically, the generator $x\otimes t_\zeta^n$ acts as $x[n]$ on the
module $V_\lambda(\zeta)$. In the physics literature \cite{BPZ}, this
action is sometimes denoted by $x_n(\zeta)$.  Equivalently, in terms of the
generating current $x(z) = \sum_{n\in \Z}x[n] z^{-n-1} $, let $v\in
V_\lambda(\zeta)$. The action of $x\otimes t_\zeta^n$ may be written as
$$
x\otimes t_\zeta^n \cdot v =\frac{1}{2\pi i} \oint_{C_\zeta}dz (z-\zeta)^n
x(z) v,
$$
where $C_\zeta$ is a contour around $\zeta$.

The central extension of $\widetilde{\g}_{(\zeta)}$ is isomorphic to
$\widehat{\g}'$, where the cocycle acts in the same way as on modules
localized at $0$:
$$
\langle x\otimes f(t_\zeta),y\otimes g(t_\zeta)\rangle =
\langle x,y\rangle \frac{1}{2\pi i} \oint_{t_\zeta=0} f'(t_\zeta)
g(t_\zeta) dt_\zeta  ,
$$
where $\langle x,y\rangle$ is the symmetric bilinear form on $\g$.
We call the centrally extended algebra with this cocycle
$\widehat{\g}'_{(\zeta)}$. Obviously, its representations are isomorphic to
those of $\widehat{\g}'$.

We also allow the point $\zeta=\infty$, and at that point we choose
the local variable to be $t_\infty = t^{-1}$.

\subsection{Fusion product of $\widehat{\g}_\zeta'$-modules}
Let $N\in \N$ and let $(\zeta_1,\ldots,\zeta_N)$ be $N$ distinct, finite
points in $\C P^1$ (for convenience we choose $\zeta_p\neq 0$).
Denote the local variable at each point by $t_p=t-{\zeta_p}$.

At each point $\zeta_p$, we localize an integrable
$\widehat{\g}_{(\zeta_p)}'$-module $V_p = V_{\mu_p}(\zeta_p)$ of level
$k$, and top component $\pi_p=\pi_{\mu_p}$.  We choose to consider only
modules with highest weights of the form $\mu_p =
a_p\omega_{\alpha_p}$, where $1\leq\alpha_p\leq r$ and $a_p\in \Z_{\geq0}$.
That is, highest weights corresponding to rectangular Young diagrams.

The completed loop algebra $\mathcal U = \oplus_p \g\otimes
\C[t_p,t_p^{-1}]\subset \g\otimes \C(t)$ acts acts on the tensor product of
these modules, $V_1\otimes\cdots\otimes V_N$ by the usual coproduct,
$$
\Delta^N_{\boldsymbol\zeta} (x\otimes f(t)) = 
\sum_{p=1}^N \left(x \otimes f(t_p+\zeta_p)\right)_{(p)},
$$
where the $p$th term in the sum above acts on the $p$th factor in the
tensor product only:
$$
x_{(p)} w_1\otimes \cdots \otimes w_N := w_1\otimes \cdots \otimes x\cdot
w_p \otimes \cdots \otimes w_N, \ x\in\mathcal U.
$$
Here, by $\C(t)$ we mean rational functions in $t$, although
we need only consider for our purposes the smaller space of rational
functions with poles at at most $\zeta_1,\ldots,\zeta_N$.

This action has a central extension, where the cocycle acts as
$$
\langle x\otimes f(t), y\otimes g(t)\rangle = \langle x,y\rangle
\sum_{p=1}^N \frac{1}{2\pi i} 
\oint_{t=\zeta_p} f'(t) g(t) dt,\ f(t), g(t)\in \C(t).
$$
Thus, the level of the action of the centrally extended, completed algebra
$\widehat{\mathcal U} = \mathcal U \oplus \C c$ is also $k$, which is
the same 
as the level of each localized module $V_i$. This action is called the
{\em fusion action} in the physics literature.  Since it differs from
the usual action on the tensor product of $\widehat{g}$-modules (which
has level $Nk$), it is denoted in \cite{FJKLM} by the symbol
$\boxtimes$ rather than the usual $\otimes$:
\begin{equation}\label{fusionproduct}
\mathbf V_{\boldsymbol\mu} (\boldsymbol\zeta) \overset{\rm def}{=}
V_{\mu_1}(\zeta_1)\boxtimes \cdots \boxtimes V_{\mu_N}(\zeta_N),\
\bmu = (\mu_1,\ldots,\mu_N).
\end{equation}

\subsection{Coinvariant spaces}
The fusion product is an integrable $\widehat{\g}'$-module of
level-$k$, thus, there is a sense in which it is completely reducible
(see Appendix I of \cite{FKLMM1} for the precise explanation and
proofs). The ``multiplicity'' of the irreducible 
$\widehat{\g}'$-module $V_\lambda(0)$ in the fusion product is given by
the Verlinde numbers \cite{TK}, which we denote by
$K_{\lambda,\bmu}^{(k)}$. If $k$ is sufficiently large (that
is, $k\geq \sum_{p}a_p$), these numbers are just the
sums of products of the usual Richardson-Littlewood coefficients. In
this paper, we only need to consider this case in order to obtain the
character formul\ae\ .

\begin{remark}
  In the case where $\alpha_p=1$ for all $p$ and $k$ is sufficiently
  large, the multiplicities are
  the usual Kostka numbers $K_{{\olambda},\mu}$ in the notation of
  \cite{Mac}, where $\mu =(a_1,\ldots,a_N)$ and ${\olambda}$ is
  a partition of length $r+1$ with $|{\olambda}| = |\mu|$,
  such that
  ${\olambda}_i-{\olambda}_{i+1}=\lambda(\alpha_i)$,
  where $\alpha_i$ are the simple roots.
\end{remark}

\newcommand{\bzeta}{{\boldsymbol\zeta}}
In complete generality, the multiplicity $K_{\lambda,\bmu}^{(k)}$ is
equal to the dimension of the coinvariant space \cite{TK,FKLMM1}
$$
\mathcal C_{\lambda, \bmu}(\bzeta) := V_{\lambda^*}(\infty) \boxtimes
\mathbf V_\bmu(\bzeta)/\langle\g\otimes \mathcal A\rangle,
$$
where the quotient is taken with respect the image of $\g\otimes
\mathcal A$ acting on the fusion product, where $\mathcal A$ is the
space of meromorphic functions with possible poles at the points $\zeta_p$
and $\infty$ (it has trivial central extension). Here, $\lambda^*$
refers to the highest weight of the dual module to $\pi_\lambda$:
$\lambda^*=-\omega_0(\lambda)$ where $\omega_0$ is the longest element
in the Weyl group.

\subsection{The coinvariant space as a quotient of principal subspaces}
The dimension of the coinvariant space $\mathcal C_{\lambda,\bmu}$
was the subject of the paper \cite{FJKLM}, where a grading was defined
on the space, compatible with the action of the current algebra. We
will use the results about this space here, and compute the graded
dimension for the special case of rectangular Young diagrams, with
sufficiently large $k$.

\begin{thm}(\cite{FJKLM} (1.6), slightly modified)
There is a surjective map
$$
u_{\lambda^*}(\infty) \otimes
\pi_1\otimes\cdots\otimes\pi_N  \to
\mathcal C_{\lambda,\boldsymbol \mu}(\bzeta),
$$
where $u_{\lambda^*}(\infty)$ is the lowest weight vector of the top
component of
the module $V_{\lambda^*}(\infty)$ with respect to the action of $\g$, 
and $\pi_p$ are the top
components of the modules $V_{\mu_p}(\zeta_p)$.
\end{thm}

Thus we can conclude that the coinvariant is a quotient of the fusion
product of principal subspaces $W_p = W_{\mu_p}(\zeta) =
U(\n_-\otimes\C[t_p^{-1}])v_p$, where $v_p$ is the highest-weight
vector of $V_p$, because
$\pi_p\subset W_p$.
The fusion product of principal subspaces is the space
$$\mathbf W_\bmu(\bzeta) = W_1\boxtimes \cdots \boxtimes W_N =
U(\n_-\otimes \C(t)) v_1\otimes\cdots\otimes v_N,$$
where we allow poles at $t=\zeta_p$.

That is, in exactly the same way as for the integrable modules, the
fusion product of principal subspaces can be decomposed as a
direct sum of principal subspaces $W_\lambda(0)$, with multiplicities
given by the Verlinde numbers $K_{\lambda, \bmu}^{(k)}$.

We can compute these multiplicities by computing the dimension of the
space of highest-weight vectors (with respect to the action of $\g$) in
the space $U(\n_-\otimes \C[t])v_1\otimes\cdots\otimes v_N$. Notice
that $x\otimes t^n$ acts on the $p$th factor by $x \zeta_p^n$.  
(Here, we do not allow poles at $\zeta_p$, because they generate vectors
in $W_p$ which are not in the top component $\pi_p$.)

\begin{remark}
The naturally graded version of the space described in the previous paragraph
is the Feigin-Loktev ``fusion product'' \cite{FL}.
\end{remark}

\subsection{Dual space of functions to the coinvariant}
Again, in this paper, we do not incorporate the level-restriction for
$k$, but we simply assume $k$ to be sufficiently large, with respect to
the collection of weights $\mu_p$: if $\mu_p = a_p \omega_{\alpha_p}$,
then the assumption is equivalent to $k\geq\sum_p a_p$. In this case, the
Verlinde number $K_{\lambda,\boldsymbol{\mu}}^{(k)}$ is equal to the
Littlewood Richardson coefficient $K_{\lambda,\boldsymbol{\mu}}$. This is
all we need in this paper to compute the characters of $W_\lambda$ for
generic $\lambda\in P^+_k$.

Consider the space of matrix elements
$C_{\lambda,\bmu}$, also known as the
space of conformal blocks:
\begin{equation}\label{matrixelement}
C_{\lambda,\bmu}=\left\{\langle u_{\lambda^*}| U(\n_-\otimes\C[t])
v_1(\zeta_1)\otimes\cdots\otimes v_N(\zeta_N)\rangle.\right\}
\end{equation}
Here, $u_{\lambda^*}$ is the lowest weight vector of
$V_{\lambda^*}(\infty)$, considered as $\widehat{\g}_{(0)}$-module
with $\widehat{\g}_{(0)}$ acting to the left. (Thus, $\n_-\otimes
\C[t^{-1}]$ acts on $u_{\lambda^*}$ trivially.)  

If $\zeta_j$ are pairwise distinct, the action of $\n_-\otimes \C[t]$
on the product of highest-weight vectors generates all of
$\pi_1\otimes\cdots\otimes
\pi_N$ (c.f. the fusion product of
\cite{FL}). The multiplicity of $v_\lambda\in\pi_\lambda$ in this
tensor product is the Littlewood Richardson coefficient $K_{\lambda,\bmu}$.

This space has a filtration by degree in $t$ inherited from the
corresponding filtration on the universal enveloping algebra.  Let
$U^{\leq n}$ be the subspace of elements in $U(\n_-\otimes\C[t])$ of
degree less than or equal to $n$ in $t$.  Let $C_{\lambda,\bmu}^{\leq
n}$ be the subspace of matrix elements of $U^{\leq n}$.
Let $C_{\lambda,\bmu}[n]= \Gr_n C_{\lambda,\bmu}$ be the
graded component of degree $n$. We define the graded coefficients
$\mathcal K_{\lambda,\bmu}(q^{-1})$ to be
\begin{equation}
\mathcal K_{\lambda,\bmu}(q^{-1}) = \sum_n q^{-n}\dim C_{\lambda\bmu}[n].
\end{equation}
We choose powers of $q^{-1}$ rather than $q$ in order to be consistent
with the grading in the last section, where we defined the degree of
$f[n]$ to be $-n$, as in \eqref{degreedef}. Therefore,
$\cK_{\lambda,\bmu}(q)$ is a polynomial in positive powers of $q$.
(Notice that this is by definition the coefficient of $\pi_\lambda$ in
the fusion product of Feigin and Loktev \cite{FL}.)

Let $\mathcal G(\boldsymbol{\zeta})_{\lambda,\boldsymbol{\mu}}$ be the
space of generating functions for matrix elements of the form
(\ref{matrixelement}). That is,
\begin{equation}\label{me}
\G = \left\{ \langle u_{\lambda^*}
|f_{\alpha_1}(x_1^{(\alpha_1)})\cdots f_{\alpha_m}
(x^{(\alpha_m)}_{m^{(\alpha_m)}}) 
v_1(\zeta_1)\otimes\cdots\otimes v_N(\zeta_N)\rangle\right\},\quad
\end{equation}
where $f_\al(x) = \sum_n f_\al[n] x^{-n-1}$
and $1\leq \al\leq r$.

Obviously, for this matrix element to be non-zero, the sum of the
$\h$-weights should be 0, that is, the matrix element should be
$\g$-invariant. 
If there are exactly $m^{(\al)}$ generating currents of the form
$f_{\al}(x_i^{(\al)})$ in the matrix element (\ref{me}), define
$\bm=(m^{(1)},\ldots,m^{(r)})^T$. Then $\bm$ is fixed by the zero-weight
condition on the matrix element. Specifically, let
$\bomega=(\omega_1,\ldots,\omega_r)^T$. Then the zero-weight condition on
$\bm$ is
\begin{equation}\label{zeroweight}
\sum_p \mu_p - \bomega^T C_r \bm  - \lambda = 0.
\end{equation}
Recall the notation $\lambda = \sum_\al l_\al \omega_\al$, and
$\mathbf l = (l_1,\ldots,l_r)^T$.
Let 
$$n_a^{(\al)}= \hbox{number of weights of the form } \mu_p=a \omega_\alpha,$$ 
and  $n^{(\al)} = \sum_a a n_a^{(\al)}$, $\bn = (n^{(1)},\ldots,n^{(r)})^T$.
Then $\sum_p
\mu_p = \sum_\al n^{(\al)}\omega_\al$. We can rewrite
(\ref{zeroweight}) more compactly as
\begin{equation}\label{mrestriction}
\bm = C_r^{-1}(\bn -\mathbf l),
\end{equation}
where $C_r$ is the Cartan matrix of $\sl_{r+1}$.

Let $g(\mathbf x)\in \G$, where $\mathbf x =
\{x_i^{(\alpha)}, i=1,\ldots,m^{(\alpha)}; \alpha=1,\ldots,r\}$.  We define
the pairing between functions $g(\mathbf x)$ and an element in
$U(\n_-\otimes\C(t))$
of the form $M (f_\alpha\otimes t_p^n)_{(p)}$, where
$M\in U(\n_-\otimes\C[t])$ and $x_{(p)}$ is an element in the algebra which 
acts on the $p$th factor only. The pairing 
is again defined inductively as in (\ref{pairing}), but the integral is
modified to
\begin{equation}\label{newpair}
(g(\mathbf x), M (f_\alpha\otimes t_p^n )_{(p)}) =
\Bigl(  \frac{1}{2 \pi i}\oint_{\mathcal C_p} g(\mathbf x)
(x_1^{(\alpha)}-\zeta_p)^n 
dx_1^{(\alpha)}, M \Bigr)
\end{equation}
where $\mathcal C_p$ is a contour around the point
$\zeta_p$, and so forth.

We now describe the zero and pole structure of the space of functions
$\G$. First we note that
$\G$ is a subspace of the dual space $\mathcal G[\bm]$ to
$U(\n_-\otimes \C[t,t^{-1}])[\bm]$, which is described in Theorem
\ref{dualtoU}:
$$
\mathcal G[\bm] =\left\{
  \left. \frac{g_1(\bx)}{\prod(x_i^{(\al)}-x_j^{(\al+1)})}\ \right| \ 
g_1(\bx)|_{x_1^{(\al)}=x_2^{(\al)}=x_1^{(\al\pm1)}}=0,\
g_1(\bx)|_{x_i^{(\al)}\leftrightarrow x_j^{(\al)}} = g_1(\bx).
\right\} .
$$

Also, recall that $f_\alpha\otimes t_p^0 = f_\alpha[0](\zeta_p)$ acts
trivially on $v_p$, unless $\mu_p=a_p\omega_\alpha$, in which case,
$(f_{\alpha}[0])^{a_p+1} $ acts trivially on $v_p$.  In addition,
$f_\alpha\otimes t_p^n$ acts trivially on $v_p$ for all $n>0$ and all
$\alpha$.

This implies, from the pairing (\ref{newpair}), that for $g(\mathbf
x)\in \G$, ${g}_1(\mathbf x)$ in equation (\ref{rational}) can have at
most a simple pole whenever $x_i^{(\alpha_p)}=\zeta_p$. There is no
pole when $x_i^{(\beta)}=\zeta_p$ if $\beta\neq \alpha_p$.
That is,
\begin{equation}\label{fracform}
g(\mathbf x) = \frac{{g}_2(\mathbf
  x)}{\prod(x_i^{(\alpha)}-x_j^{(\alpha+1)}) \prod_{p}
    (x_a^{(\alpha_p)}-\zeta_p)}\in\G,
\end{equation}
where the function ${g}_2(\mathbf x)$ satisfies
\begin{equation}\label{zetavanishing}
{g}_2(\mathbf x)|_{
x_1^{(\alpha_p)} = \cdots = x_{a_p+1}^{(\alpha_p)} = \zeta_p}=0 
\ , \qquad \forall p \ .
\end{equation}
(Recall that we assume $\zeta_p\neq 0$, so that there is no pole at
$x_i^{(\alpha)}=0$.) Thus, ${g}_2(\mathbf x)$ is a polynomial in
$x_i^{(\alpha)}$.

Finally, the currents in $U(\n_-\otimes \C[t])$ may act to the
left, on $u_{\lambda^*}$ sitting at infinity. The pairing at infinity is
\begin{eqnarray*}
( g(\mathbf x), (f_\alpha \otimes t^n)_{(\infty)} M ) &=&
\Bigl( \frac{1}{2\pi i}  \oint_{C_\infty} g(\mathbf x)
(x_1^{(\alpha)})^{n} dx_1^{(\alpha)}, M \Bigr) \\
& =&  \Bigl( \frac{1}{2\pi i}\ \oint_{C_0}
(x_1^{(\alpha)})^{-n-2}g((x_1^{(\alpha)})^{-1},x_2^{(\alpha)},\ldots)
d x_1^{(\alpha)}, M \Bigr)
\end{eqnarray*}
(the contour around infinity is clockwise).
Since $f_\alpha[n]$ acts trivially at $\infty$ if $n\leq 0$, this
integral should be zero for $n\leq 0$ if $g(\mathbf x)\in\G$.
This shows that
\begin{equation}\label{degree}
\deg_{x_i^{(\alpha)}}g(\bx) \leq -2 \qquad \hbox{for all $i,\alpha$}.
\end{equation}

In summary, we have that, for $k$ sufficiently large,
\begin{thm}
  The dual space $\G$ to the space coinvariants $\mathcal
  C_{\lambda,\boldsymbol{\mu}}(\bzeta)$, with respect to the pairing
  (\ref{newpair}), is the space of functions in the
  variables $\mathbf x = \{x_i^{(\alpha)}\ |\ \alpha=1,\ldots,r;
  i=1,\ldots,m^{(\alpha)}\} $, where $m^{(\alpha)}$ is determined by
  (\ref{mrestriction}),
  of the form (\ref{fracform}), where ${g}_2 (\mathbf x)$ is a
  polynomial, symmetric with respect to exchange of variables with the
  same superscript $(\alpha)$, satisfying the Serre relation
  (\ref{serre1}) and the vanishing condition
  (\ref{zetavanishing}), with the degree of ${g}(\mathbf x)$ in each
  variable less than or equal to $-2$.
\end{thm}

In the next section, we compute the character of this space.

\subsection{Filtration of the dual space}
The space $\mathcal G_{\lambda,\boldsymbol\mu}(\boldsymbol\zeta)$ is
filtered by homogeneous (total) degree in $x_i^{(\alpha)}$. Let
$\G[n]$ be the graded component. This space is dual to the space
$\mathcal C_{\lambda,\bmu}[n+|\bm|]$ (because the definition of the
pairing involves taking the residue).  We normalize the degree of the
cyclic vector to be 0.
Therefore we have
\begin{equation}
\ch_q C_{\lambda,\bmu}=\ch_q\mathcal C_{\lambda,\bmu} = q^{-|\bm|}\ch_q\G = 
\sum_{n} q^{-n-|\bm|} \G[n] = \mathcal
K_{\lambda,\bmu}(q^{-1}). 
\end{equation}

We use the same filtration argument as in Section 3. That is, consider
the lexicographic ordering on $r$-tuples of partitions
$\boldsymbol\nu$, where $\nu^{(\alpha)}$ is a partition of
$m^{(\alpha)}$. (Since $k$ plays no role in the filtration argument
except in limiting the types of partitions allowed in the filtration,
there is no difference in the zero and pole structure related to
$k$). We act with the evaluation maps $\phi_{\boldsymbol{\nu}}$ on the
space $\mathcal G_{\lambda,\boldsymbol\mu}(\boldsymbol\zeta)$ and
consider the image in the space $\mathcal H[\bm]$ of functions in the
variables
$$
\{y_{a,i}^{(\alpha)}\ |\ a\geq 1, i=1,\ldots,m_a^{(\alpha)}, m_a^{(\alpha)}
= {\rm Card}(\{\nu^{(\alpha)}_i=a\})\}$$
 of 
the subspaces $\Gamma_{\boldsymbol\nu} =
\cap_{\boldsymbol\nu'>\boldsymbol\nu}{\rm Ker}\phi_{\bnu'}$. 
We take the associated graded
space, and compute the character of the graded components 
$\Gamma_\bnu/\Gamma'_\bnu$, where $\Gamma'_\bnu=\cap_{\bnu'\geq\bnu}
{\rm Ker}\phi_{\bnu'}.$ Define $\mathcal H_{\bnu}$ to be the image of
the induced map ${\overline{\varphi}_\nu}: \Gamma_\bnu/\Gamma'_\bnu$.

The results are as follows.
\begin{lemma}\label{Gzeros}
Let $g(\mathbf x)\in\G$. Then 
$$\phi_{\boldsymbol\nu}(g(\mathbf x)) = \prod_{\alpha; (a,i)<(a',i')}
(y_{a,i}^{(\alpha)}-y_{a',i'}^{(\alpha)})^{2A_{a,a'}}
h_1(\mathbf y).$$
\end{lemma}
\begin{proof}
This follows from Lemma \ref{lemintzeros}. The only
difference in the two situations is that the partitions are only
restricted by $m^{(\alpha)}$, not $k$.
\end{proof}

The next Lemma gives the pole structure due to the nontrivial commutation
relations together with the Serre relations. Its proof is identical to Lemma
\ref{lemserpoles}. 
\begin{lemma}\label{serrezeros}
Let $h_1(\mathbf y)$ be defined as in Lemma \ref{Gzeros}. Then
$$
h_1(\mathbf y) = \prod_{\alpha = 1}^{r-1} \prod_{a,a',i,i'}
(y_{a,i}^{(\alpha)} - y_{a',i'}^{(\alpha+1)})^{-A_{a,a'}} h_2(\mathbf y),
$$
where $h_2(\mathbf y)$ is regular when $y_{a,i}^{(\alpha)} =
y_{a',i'}^{(\alpha+1)}$.
\end{lemma}

The following Lemma is a slight modification of Lemma \ref{lemirrpoles}.
\begin{lemma}
Let $h_2(\mathbf y)$ be as in Lemma \ref{serrezeros}. Then
$h_2(\mathbf y)$ as a pole of order at most 
$\min(a,a_p)$ whenever $y_{a,i}^{(\alpha_p)} = \zeta_p$.
\end{lemma}

Thus, we have
\begin{equation}\label{hofy}
\phi_{\boldsymbol\nu}(g(\mathbf x)) = h(\mathbf y) 
= \frac{\prod (y_{a,i}^{(\alpha)}-y_{a',i'}^{(\alpha)})^{2A_{a,a'}}}
{\prod (y_{a,i}^{(\alpha)}-y_{a',i'}^{(\alpha+1)})^{A_{a,a'}}}
\prod_p
(y_{a,i}^{(\alpha_p)}-\zeta_p)^{-A_{a,a_p}} h_3(\by),
\end{equation}
where $h_3(\mathbf y)$ is a polynomial in the variables
$\{y_{a,i}^{(\alpha)}\ |\ a\geq 1,
i=1,\ldots,m_a^{(\alpha)},\alpha=1,\ldots,r\}$, with $\sum_a a
m_a^{(\alpha)} = m^{(\alpha)}$, symmetric under the exchange of
variables $y_{a,i}^{(\alpha)}\leftrightarrow
y_{a,i'}^{(\alpha)}$. Here, $m_a^{(\alpha)}$ is the number of parts
of length $a$ in the partition $\nu^{(\alpha)}$.

\begin{remark}
It is important to note that, since we are only interested in the
character of the space of functions of the form (\ref{hofy}), we can
now set all $\zeta_p=0$ in the space of polynomials without changing
the character of the space.
\end{remark}

There is a further restriction on $h(\by)$ coming from the degree
restriction (\ref{degree}) on $g(\bx)$. (This ensures
that the space of coinvariants is
finite-dimensional.) The
evaluation map is degree preserving, which implies
that $$\deg_{y_{a,i}^{(\alpha)}}h(\by)\leq -2a.$$ This gives
the following restriction on the degree of $h_3(\mathbf y)$:
\begin{lemma}\label{deglemma}
Let $h_3(\mathbf y)$ be as in equation (\ref{hofy}). Then $h_3(\mathbf
y)$ is a polynomial in the variables $\{y_{a,i}^{(\alpha)}\}$, with 
$$
0\leq deg_{y_{a,i}^{(\alpha)}} h_3(\mathbf y) \leq
-\sum_{b,\beta}(C_r)_{\alpha,\beta}A_{a,b}m_b^{(\beta)} + \sum_b
A_{a,b}n^{(\alpha)}_b,
$$
where $n_a^{(\alpha)}$ is the number of $\g$-modules with highest
weight $a \omega_\alpha$.
\end{lemma}
\renewcommand{\H}{{\mathcal H_{\boldsymbol\nu}}}

The injectivity of the induced map $\overline{\varphi}_\nu:
\Gamma_\nu/\Gamma'_\nu \to \mathcal H_\bnu$ 
follows from the injectivity argument of Lemma
\ref{injective}. 

We do not show surjectivity. Instead, we compute the graded character
of the coinvariant using the above space of functions, evaluate it at
$q=1$, and show that it is equal to the desired multiplicity given by
the Littlewood-Richardson rule, by comparing with the known result
\cite{KSS} for generalized Kostka polynomials.  

The argument is as follows. The injectivity of the map
$\overline{\varphi}_\bnu$, which is a degree preserving map, implies
that $$\dim \G[n] \leq \sum_\bnu\dim \H[n],$$
where by $[n]$ means the
graded component with respect to the homogeneous grading in the
variables $\by$.  We will show that $\dim\G = \sum_\bnu \dim\H$, by
computing the $q$-character of $\H$, and showing that $\dim\H=
K_{\lambda,\bmu}$, which is the dimension of the space of
coinvariants. This proves the surjectivity of the evaluation map
$\overline{\varphi}_\bnu$, and also gives the $q$-character of $\G$.

Define the character of the space $\H$
to be 
$$
\ch_q\H =  \sum_n q^{-n} \H [n].
$$
This character can be computed by
setting $\zeta_p\to 0$ for all $p$. Recall that we must multiply by
$q^{-|\bm|}$ to obtain the character of the coinvariant. We use the 
Gaussian polynomial,
$$\qbin{m+n}{m}_q = \frac{(q)_{m+n}}{(q)_{m}(q)_{n}},\quad
m,n\in\Z_{\geq 0}.$$ 

\begin{lemma}
Let $\H$ be the space of functions of the
form (\ref{hofy}) with degree restrictions \eqref{deglemma}, and
$\zeta_p=0$. Then
$$
q^{-|\bm|}\ch_q\H = q^{Q(\mathbf m,\mathbf
  n)}\prod_{a,\alpha}
\qbin{P_a^{(\alpha)}+m_a^{(\alpha)}}{m_a^{(\alpha)}}_q,
$$
where 
$$
Q(\mathbf m,\mathbf n) = 
\frac{1}{2}\sum_{a,b,\alpha,\beta} m_a^{(\alpha)}
(C_r)_{\alpha,\beta}A_{a,b}m_b^{(\beta)} - \sum_{a,b,\alpha}m_a^{(\alpha)}
A_{a,b} n_b^{(\alpha)}
$$
and 
$$
P_a^{(\alpha)} = \sum_{b} A_{a,b} n_b^{(\alpha)} -
\sum_{\beta,b} (C_r)_{\alpha,\beta}A_{a,b} m_b^{(\beta)}.
$$
Here, $m_a^{(\alpha)}$ is the number of parts of $\nu^{(\alpha)}$ of
length $a$, and $n_a^{(\alpha)}$ is the number of $\g$-modules of
highest weight $a\omega_\alpha$.
\end{lemma}

Since the evaluation maps $\phi_{\boldsymbol\nu}$ are degree
preserving, we can conclude that
$$
\ch_q \G \leq \sum_{\boldsymbol\nu} ch_q \H,
$$
where by the inequality, we mean the inequality in the coefficient of
each power of $q$.

Recall the identity
$$\qbin{m+n}{m}_q = q^{m n}\qbin{m+n}{m}_{\frac{1}{q}}.$$
We can now conclude that we have an equality.
\begin{thm}\label{generalizedkostka}
The graded character of the space of conformal blocks
$C_{\lambda,\bmu}$
is $\mathcal K_{\lambda,\bmu}(q^{-1})$, where
\begin{equation}\label{Kostka}
\mathcal{K}_{\lambda,\boldsymbol\mu}(q) = 
\sum_{\vbm} q^{\frac{1}{2} \vbm^T C_r\otimes A_k \vbm}
  \prod\qbin{P_a^{(\alpha)}+m_a^{(\alpha)}}{m_a^{(\alpha)}}_q
\end{equation}
where $\vbm$ is a vector with entries $m_a^{(\al)}$ restricted by
(\ref{mrestriction}), namely 
$\bm = C^{-1}_r (\bn-\bl)$, and 
$$
\overset{\to}{\mathbf P} = (\id\otimes A_k)\vbn - (C_r\otimes A_k)\vbm.
$$
\end{thm}
\begin{proof}
A direct comparison of the fermionic formula on the right
hand side of (\ref{Kostka}) with equation (2.6) of \cite{KSS} shows that 
\begin{equation}
\label{kcomp}
\cK_{\lambda,\bmu} (q)=
K_{\overline\lambda^t,R^t} (q) \ ,
\end{equation}
in the notation of \cite{KSS} (where $K_{\lambda,R}(q)$ is the
co-charge Kostka polynomial). Here, $\overline\lambda$ is the Young
diagram obtained from the weight
$\lambda$ by adjoining to the corresponding Young diagram of $\lambda$
columns of length $r+1$, so that the equality
$|\overline\lambda| = |R|$ is satisfied (the Kostka polynomial is zero
unless $|R|-|\lambda|\equiv 0 \mod(r+1)$,  as a consequence of the
restriction on the summation over $m_a^{\alpha}$, see
part (4) of Lemma \ref{kostkaprops} below). The sequence
$R=(R_1,\ldots,R_N)$, with $R_p = (a_p)^{\alpha_p}$, is the sequence of
rectangular Young diagrams corresponding to the weights $\mu_p$. 

We use a duality theorem for generalized Kostka polynomials \cite{Kir}
\begin{equation}
\label{ktr}
K_{\lambda^t;R^t} (q)=
q^{n(R)} K_{\lambda,R} (q^{-1}) \ ,
\end{equation}
where $n(R) =\sum_{1\leq p<p'\leq N}
\min(\alpha_p,\alpha_{p'})\min(a_p,a_{p'})$.  Then using the fact that
$$K_{\lambda,R}(1) = \dim Hom_{\g}(\pi_\lambda,
\pi_{\mu_1}\otimes\cdots\otimes\pi_{\mu_N})$$
(where $\g=\sl_{r+1}$ or
$\gl_{r+1}$) is the dimension of the space of conformal blocks $
C_{\lambda\boldsymbol\mu}$, we conclude the equality of $q$-dimensions
in the Theorem holds.
\end{proof}

\subsection{A remark about the structure of $\cK_{\lambda,\bmu}(q)$}
\label{kostkastructure}
In this paper, since we are concerned with representations of
$\sl_{r+1}$, we have labeled the representations with highest weight
$\lambda$ with respect to the $\sl_{r+1}$ weights, $\lambda =
l_1\omega_1 + \cdots + l_r\omega_r$, and similarly for the weights
$\mu_p = a_p \omega_{\alpha_p}$ with $\alpha_p\leq r$.

Define $\mathcal S_N^{(k)}$ to be the set of all unordered $N$-tuples of
$\sl_{r+1}$ dominant weights of the form $\mu_p = a_p \omega_{\alpha_p}$,
with $\sum_p a_p\leq k$.  Let $P(r,k)$ be the set of partitions of
length at most $r$ and width at most $k$. Define $\nu:\mathcal
S_N^{(k)} \to P(r,k)$ to be the ``horizontal concatenation'' map:
$$(\nu(\bmu))_\beta = \sum_{\alpha=\beta}^r\bn^{(\alpha)},
\quad 1\leq \beta\leq r, \bmu\in\mathcal S_N^{(k)}.$$
Note that this map is surjective but in general not injective.

Let $\mathcal S_r$ be the subset of $S_r^{(k)}$
consisting of precisely $r$ weights of the form $\mu_p= a_p \omega_p$
(again with $\sum a_p \leq k$). That is, $n^{(\alpha)}=a_\alpha$.
 Then $\nu$ is now a natural isomorphism,
$\nu: \mathcal S_r\overset{\sim}{\to} P(r,k)$.
The inverse map is $\nu^{-1}(\mu) = \bmu=(\mu_1,\ldots,\mu_r)$, with
$\mu_p = (\mu_p-\mu_{p+1})\omega_p$, with $\mu_{r+1} = 0$ by
definition. 

In this paper we need to consider only the cases where $\bmu\in
\mathcal S_r$ and $\lambda\in P(r,k).$ In this special case, we
have the following properties of the Kostka polynomial.
\begin{lemma}\label{kostkaprops}
Let $\bmu\in \mathcal S_r$ and $\lambda\in P(r,k)$.
Then the following statements are true for the Kostka polynomial of
equation \eqref{Kostka}:
\begin{enumerate}
\item $\cK_{\lambda,\bmu}(q)
 = 1$ if $\nu(\bmu)=\lambda$;
\item $\cK_{\lambda,\bmu}(q)=0$ if $\lambda_1 > \nu(\bmu)_1$;
\item  $\cK_{\lambda,\bmu}(q)=0$ if $\lambda_1=\nu(\bmu)_1$ and
 $\lambda_s > \nu(\bmu)_s$ where $s$ is the smallest integer such that
$\lambda_s \neq \nu(\bmu)_s$;
\item $\cK_{\lambda,\bmu}(q)=0$ if $\frac{1}{r+1}(|\nu(\bmu)| -
  |\lambda|)\notin \Z_{\geq 0}$. 
\end{enumerate}
\end{lemma}

Let $\KK(q)$ be the matrix with entries $(\KK(q))_{\lambda,\nu(\bmu)}
= \cK_{\lambda,\bmu} (q)$ with $\bmu\in \mathcal S_r$ and
$\lambda\in P(r,k)$. The Lemma implies in particular that $\KK(q)$ is
upper unitriangular with respect to the ordering on partitions which
looks like the lexicographic ordering on partitions, applied to
partitions which are not of the same size: $\lambda<\mu$ if $\lambda_i
= \mu_i$ for all $i<s\leq r$, and $\lambda_s < \mu_s$.

\begin{proof}
\begin{enumerate}
\item The constraint $C_r\bm = \bn-\bl$ means,
when $\bn=\bl$ (i.e. $\lambda =\nu(\bmu)$), that $m^{(\alpha)}=0$ for
all $\alpha$, hence only the term with $m_a^{(\alpha)}=0$ contributes
to the sum.

\item
\label{item2}
Suppose that $\lambda_1 > \nu(\bmu)_1$, which implies that
$\sum_{\al=1}^r (n^{(\al)}-l_\al)<0$. However, the constraint implies that
\begin{equation}
\sum_{\al=1}^r (n^{(\al)}-l_\al)= \sum_{\alpha=1}^{r} (C_{r} \bm)_\alpha=
m^{(1)}+m^{(r)}
\end{equation}
and since $\ma \geq 0$ for non-zero Kostka polynomials, this
gives the desired result. 

\item
Arguments similar to the proof of item \eqref{item2} show that
$\lambda_1 = \nu(\bmu)_1$ implies $m^{(1)}=m^{(r)} =0$, and also
that $s<r$. Note that from $\lambda_\alpha = \nu(\bmu)_\alpha$ for $\alpha\leq s-1$ it follows
that $n^{(\alpha)}=l_\alpha$ for $\alpha\leq s-2$. 
>From the constraint, we now obtain the relations
$$
\sum_{\alpha=1}^{t} (C_r \bm)_\alpha = m^{(1)} + m^{(t)} - m^{(t+1)}=m^{(t)} - m^{(t+1)}
= 0 \qquad 1 \leq t \leq s-2 \ ,
$$
which imply that $m^{(t)} = 0$ for $1\leq t \leq s-1$, in order for the
Kostka polynomials to be non-zero. From the assumption that
$\lambda_s > \nu(\bmu)_s$, we obtain $l_{s-1} < n^{(s-1)}$. Thus, we find that
$$
\sum_{\alpha=1}^{s-1} (C_r \bm)_\alpha = m^{(1)}+ m^{(s-1)} - m^{(s)} = 
-m^{(s)} = n^{(s-1)}-l_{s-1} > 0 \ ,
$$
which implies that the Kostka polynomial indeed vanishes, because
$m^{(s)}<0$.

\item This comes from the fact that $m^{(r)}\in\Z_{\geq 0}$, and
\begin{equation}\label{mr}
m^{(r)} =
\sum_{\alpha=1}^r (C_r^{-1})_{r,\alpha}(n^{(\alpha)}-l_\alpha) =
\frac{1}{r+1}\sum_{\alpha=1}^r \alpha(n^{(\alpha)}-l_\alpha) =
\frac{1}{r+1}(|\nu(\bmu)|-|\lambda|).
\end{equation}
\end{enumerate}
\end{proof}

We can also make contact with the usual combinatorial notation for
Kostka polynomials, which are labeled by Young diagrams, that is,
$\gl_{r+1}$ representations.
Let $\olambda$ be the partition of length at most $r+1$, obtained from
$\lambda$ by defining
$$
\olambda_\beta = m^{(r)} + \lambda_\beta,
 \ 1\leq \beta\leq r+1,
 $$
 where $\bm = C_r^{-1}(\bn-\bl)$.  Let $\omu=\nu(\bmu)$. Then
 equation \eqref{mr} implies $|\olambda| = |\omu|$, which is the usual
 condition in the Kostka polynomial labeled by $\gl_{r+1}$-weights.
The partition $\olambda$ can be pictured as that obtained by adding
$m^{(r)}$ columns of length $r+1$ to the left of the Young diagram
corresponding to $\lambda$.

The Kostka polynomial is defined for any $r$. If we choose to fix
$|\omu|=m$, and choose $r$ sufficiently large ($r\geq m$), then
$n^{(\al)}=l_\al=0$ if $\al>m$. 
We have the following generalization of the  triangularity
property for Kostka polynomials: 
\begin{lemma}\label{triangular}
The generalized Kostka polynomial $\cK_{\lambda,\bmu}(q)=0$ unless
$\olambda \unlhd \omu=\nu(\bmu)$ according to the dominance ordering
on partitions.
\end{lemma}
\begin{proof}
The dominance ordering on partitions is
$$
\sum_{\alpha=1}^\beta \olambda_\alpha \leq \sum_{\alpha=1}^\beta
\omu_\alpha,\quad \hbox{for all }\beta\in 1,\ldots,r+1.
$$
Recast in terms of the variables $\bn$ and $\bl$ this means that
\begin{equation}\label{step2}
A(\bn - \bl)_\beta - \beta l_{r+1} = A(\bn - \bl)_\beta - \beta m^{(r)}
\geq 0\quad \hbox{for all }\beta.
\end{equation}
For $\beta=r+1$ the equality holds due to the condition $|\olambda|=|\omu|$,
so we need only consider $\beta\leq r$. Using the fact that
$$
m^{(r)} = 
\frac{1}{r+1}\sum_{\alpha=1}^r \alpha(n^{(\alpha)} - l_\alpha),
$$
the equation (\ref{step2}) becomes
$$
\sum_{\alpha=1}^r (A_{\beta\alpha} -
\frac{\beta\alpha}{r+1})(n^{(\alpha)}-l_\alpha) = (C_r^{-1}(\bn -
\bl))_\beta = m^{(\beta)}.
$$
Since $m^{(\beta)}\geq 0$ in the summation in $\cK_{\lambda,\bmu}(q)$,
this proves the Lemma.
\end{proof}

Note also that if $m^{(\beta)} = 0$ for all $\beta$, then
$\olambda = \omu$. In that special case, $\cK_{\lambda,\bmu}(q)=1$.

To tie in with the usual notion of the unitriangularity of the Kostka
matrix, let $\mathcal S_r[m]\sim P(r,k)[m]$ be the subsets of
(multi-) partitions of $m$, and fix $\max(k,r)\geq m$. The number of
elements of both sets is the number of partitions of $m$. Let
$\olambda\vdash m$.  The last Lemma implies that the square matrix
$\KK(q)$, with entries indexed lexicographically by the partitions
$\nu(\bmu)$ with $\bmu\in S_r[m]$ and $\olambda$ is upper
unitriangular. That is, define
$$
(\KK(q))_{\olambda,\omu} = \cK_{\lambda,\bmu},\quad \bmu\in\mathcal
S_r[m],\ \omu=\nu(\bmu),\ \max(k,r)\geq |\olambda|=|\omu|.
$$
Then $\KK(q)_{\olambda,\omu} = 0$ if $\olambda \rhd \omu$, and it is equal
to $1$ if $\olambda = \omu$.

In the case in which we are interested, in which $r$ is fixed and may
be smaller than $|\omu|$, we take the subset of the elements of this
matrix which have the length of $\omu$ to be at most $r$, and the
length of $\olambda$ to be at most $r+1$.

\section{Characters for arbitrary highest-weight $\widehat{\sl}_{r+1}$-modules}

Let $\lambda\in P_k^+$ and let $V_\lambda$ be the highest-weight
$\widehat{\sl}_{r+1}$-module of level $k$. We are interested in computing
a fermionic formula for the character of this space, for arbitrary
$\lambda$, similar in form to the one found in Section \ref{feig-sto}.
 
We compute this character in several steps. First, we compute the
character of the fusion product of several principal subspaces
corresponding to rectangular highest weights $\mu_p$. We then use a
Weyl translation to find the character of the fusion product
of integrable modules corresponding to the same highest weights.

At this point, we choose a very particular set of $r$ rectangular
highest weights, of the form $\mu_p = a_p\omega_p$ with $p=1,\ldots,r$.
We use the decomposition of the
fusion product into the graded sum over irreducible highest-weight
modules, with coefficients given by the generalized Kostka
polynomials. This means that the character of the fusion product is
the sum over characters of irreducible modules, with coefficients
given by the Kostka polynomial.

This relation between the characters is invertible, so we use it to
write the character of the irreducible module in terms of a finite sum
over characters of particular fusion products. The coefficients in the
sum are polynomials in $q^{-1}$ whose coefficients are not necessarily
positive, since they are given by the entries of the inverse of the
matrix of generalized Kostka polynomials in $q^{-1}$.

\subsection{Character of the fusion product of principal
subspaces}

Consider the fusion product of principal subspaces:
$$\bW_{\bmu}(\boldsymbol{\zeta})=W_1 (\zeta_1)
\boxtimes \cdots \boxtimes W_N (\zeta_N)= U(\n_-\otimes \C(t))
v_1\otimes \cdots \otimes v_N,
$$
where we allow singularities at $t=\zeta_p$.  Here, $v_p$ is the
highest-weight vector of $V_{\mu_p}(\zeta_p)$, the module of
level-$k$, with highest weight of the form $\mu_p=a_p \omega_{\alpha_p}$,
localized at $\zeta_p$. 

We choose $k$ sufficiently large -- that is, $k\geq \sum_p a_p$, so
that the level-restriction in the decomposition coefficients does not
play a role.

Note once more that the algebra $U(\n_-\otimes\C(t))$ is filtered by
degree in $t$, and that, defining the cyclic vector $\otimes v_p$ to
have degree 0, the fusion product $\bW_\bmu(\bzeta)$ inherits this
filtration. Hence, we can define the $q$-character of
$\bW_\bmu(\bzeta)$ as the Hilbert series of the associated graded
space -- it is a Laurent series in $q$, which we can compute for
sufficiently simple $\mu_p$.

As an $\n_-\otimes\C[t,t^{-1}]$-module, $\bW_\bmu(\bzeta)$ decomposes
as a direct sum of principal subspaces $W_\lambda(0)$, with graded
coefficients which are equal to the generalized Kostka polynomials in
the previous section.  This follows from the fact that $W_\lambda(0)$
is generated by the action of $\n_-\otimes\C[t,t^{-1}]$ on the highest
weight vector of $V_\lambda(0)$, and in the previous we computed the
graded space of multiplicities of these highest-weight vectors in the
fusion product of integrable modules to be generalized Kostka
polynomials. 

Thus, we can see that
\begin{equation}\label{fdecomp}
\ch_q\bW_\bmu(\bzeta) = \sum_{\lambda}
 \cK_{\lambda,\bmu}(q^{-1})\ch_q
W_\lambda \ .
\end{equation}
Note that the sum over $\lambda$ is finite, because
$\cK_{\lambda,\bmu} (q)= 0$ when $\lambda_1>(\nu(\bmu))_1$.

In this subsection we will compute the character of the fusion
$\bW_\bmu(\bzeta)$, by characterizing the dual space of
$\n_-\otimes\C(t)$ acting on the cyclic vector $\otimes v_p$.

The dual space is the space of generating functions for matrix elements
of the form
$$
\left\{\langle w | U(\n_-\otimes\C(t)) v_1\otimes\cdots\otimes
  v_N\rangle,\ | \ w\in W_{\lambda^*}(\infty), \lambda\in P_k^{+}\right\}.
$$
Thus, the dual space $\mathcal F_\bmu(\bzeta)$ is the space of
functions in the variables $x_i^{(\alpha)}$ (with $1\leq \alpha\leq r$
and $1\leq i\leq m^{(\alpha)}$), with pairing defined in the same way
as in equation (\ref{newpair}). Thus it is the space of functions with
possible simple poles at $x_i^{(\alpha_p)}=\zeta_p$ and
$x_i^{(\alpha)} = x_j^{(\alpha\pm1)}$, such that the polynomial
$f(\bx)$ defined by
\begin{equation}
\label{rationlfx}
F(\bx) = \frac{f(\bx)}
{\prod_{p,i} (\xap_i-\zeta_p) \prod_{\alpha=1}^{r-1}\prod_{j,k}
(\xa_j-\xaa_k)} \ \in \cF_\bmu(\bzeta)
\end{equation}
is symmetric under the exchange $x_i^{(\alpha)}\leftrightarrow
x_j^{(\alpha)}$. In addition, it vanishes due to the Serre relation
whenever
$$
x_1^{(\alpha)}=x_2^{(\alpha)}=x_1^{(\alpha\pm1)}.
$$

There is no degree restriction on $f(\bx)$, since we allow for poles
at infinity in $U(\n_-\otimes \C(t))$, as well as at $t=\zeta_p$. We
do not allow for zeros at $t=\zeta_p$, so the pole structure at
$t=\zeta_p$ is as before. Moreover we have, as in the calculation of
the coinvariant, the condition that $f(\bx)$ vanishes whenever
\begin{equation}\label{irrepcx}
\xap_1=\cdots=\xap_{a_p+1} = \zeta_p \ , \quad p=1,\ldots,N \ .
\end{equation}

Finally, it is possible now to have currents $f_\alpha(z)^{k+1}$
acting non-trivially on the tensor product of highest-weight vectors.
Since $W_\lambda(0)$ is a subspace of an integrable module, where such
currents act trivially, the dual space is in the subspace which
couples trivially to such currents. That is, we must impose the
integrability condition, that $f(\bx)$ vanishes whenever
\begin{equation}
\xa_1=\cdots=\xa_{k+1} \label{intcx} \\
\end{equation}

These conditions characterize the space $\cF_\bmu(\bzeta)$. In order
to compute the character of the $\h$-graded component
$\cF_\bmu(\bzeta)[\bm]$, we introduce the same filtration as in
Section \ref{fildual}. That is, let $\bnu$ be a multi-partition
consisting of $r$ partitions, where $\nu^{(\alpha)}\vdash
m^{(\alpha)}$, (we denote this as $\bnu\vdash\bm$). We order
multi-partitions lexicographically, and introduce the evaluation maps
$\varphi_\bnu$ as in Section \ref{fildual}. The evaluation maps act on
the space $\cF_\bmu(\bzeta)$. Let $\Gamma_\bnu = \cap_{\bnu'>\bnu}\ker
\varphi_{\bnu'}$ etc., where the kernel now refers to that of the
evaluation map acting on $\cF_\bmu(\bzeta)$. Define the graded
components $\Gr_\bnu = \Gamma_\bnu/\Gamma_\bnu'$.

We compute the image of the induced map
$\overline{\varphi}_\bnu: \Gr_\bnu\to \cH_\bnu$. Here, $\cH_\bnu$ is
the space of rational functions in the variables
$$
\by = \left\{ y_{a,i}^{\alpha} \ |\ 1\leq \alpha\leq r, 1\leq i\leq
  m_a^{(\alpha)}, 1\leq a\leq k \right\},
$$
where $m_a^{(\alpha)}$ is the number of rows of length $a$ in
$\nu^{(\alpha)}$, with possible poles at $\ya_{a,i}=\yaa_{a',i'}$ and at
$\yap_{a,i} = \zeta_p$.

\begin{defn}
Let $\widetilde{\cH}_\bnu \subset \cH_\bnu$ be the
subspace of functions spanned by functions of the form
\begin{equation}
\label{h1x}
H (\by) = H_\bnu (\by) h(\by) \ ,
\end{equation}
where $h (\by)$ is a polynomial, symmetric
under the exchange of variables with the same values of $\alpha$ and
$a$, and
\begin{equation}
\label{h2x}
H_\bnu (\by) =
\prod_{\substack{\alpha=1,\ldots,r\\(a,i)>(a',i')}}
(\ya_{a,i}-\ya_{a',i'})^{2 A_{a,a'}}
\prod_{\substack{\alpha=1,\ldots,r-1\\(a,i);(a',i')}}
(\ya_{a,i}-\yaa_{a',i'})^{-A_{a,a'}}
\prod_{p,(a,i)}(\yap_{a,i}-\zeta_p)^{-A_{a,a_p}} \ .
\end{equation}
\end{defn}

By using almost identical arguments to those in Section 
\ref{fildual}, we conclude that
\begin{thm}
\label{thmigvsx}
The induced map
\begin{equation}
\label{emapgvsx}
\overline{\varphi}_{\bnu} : \Gr_{\bnu} \Gamma \rightarrow \widetilde{\cH}_\bnu
\end{equation}
is an isomorphism of graded vector spaces.
\end{thm}

Therefore we have that
$$
\ch_q\cF_\bmu(\bzeta) = \sum_{\bm}\sum_{\bnu \vdash\bm} 
\ch_q  \widetilde{\cH}_\bnu.
$$
To compute the character of $ \widetilde{\cH}_\bnu$ we can set
$\zeta_p=0$ in $H_\bnu(\by)$, as it does not change the character.
Also recall that $\ch_q \bW_\bmu(\bzeta)[\bm] = q^{|\bm|}\ch_q
\cF_\bmu(\bzeta)$. Thus we have
\begin{equation}
\label{charfp}
\ch_q \bW_{\bmu}(\bzeta) =
\sum_{\vbm \in \Z_{\geq 0}^{r\times k}}
\frac{q^{\frac{1}{2} \vbm^T (C_r \otimes A) \vbm 
- \vbm^T (\id\otimes A) \vbn}}{(q)_{\vbm}} 
e^{\bomega^T\cdot\bn-\bomega^T C_{r} \bm} \ .
\end{equation}
Recall that $\bn = (n^{(1)},\ldots,n^{(r)})^T$, with
$n^{(\alpha)} =\sum_{a\geq 0} a n^{(\alpha)}_a$, where $\na_a$ is the
number of highest weights of the form $\mu_p=a \omega_\alpha$.

In order to calculate the character for general principal
subspaces of $\widehat{\sl}_{r+1}$, we can restrict ourselves
to sequences of $r$ partitions of the form $\mu_p = a_p \omega_p$, with
$p=1,\ldots,r$. 

The results of section \ref{kostkastructure} show that the matrix $\KK (q)$
with elements $(\KK(q))_{\lambda,\nu(\bmu)}= \cK_{\lambda,\bmu} (q)$
is invertible, so we can invert the relation \eqref{fdecomp} and
conclude that the character of the principal subspace of a general
highest weight is given by
\begin{equation}
\label{wdecomp}
\ch_q  W_\lambda = \sum_{\bmu}
(\KK^{-1} (q^{-1}))_{\nu(\bmu),\lambda} \
\ch_q \bW_{\bmu}(\bzeta) \ ,
\end{equation}
where the finite sum is over sequences of partitions of the form
$\bmu = (n^{(1)} \omega_1,\ldots,n^{(r)}\omega_r)$, i.e.\
sequences of rectangular weights, such that $\nu(\bmu)\leq \lambda$
(in the sense of Lemma \ref{kostkaprops}). 

\subsection{Characters for general highest-weight modules of
$\widehat{\mathfrak{sl}}_{r+1}$}

We can now use the results of section \ref{awt}, to obtain the
character formul\ae\ for the Weyl translated principal subspaces and,
in particular, the characters of general integrable irreducible
representations of $\widehat{\mathfrak{sl}}_{r+1}$.

Let us denote the limit of $N\rightarrow\infty$ of
$T^{\bN} \ch_q \bV_{\bmu}(\bzeta)$ (where $\bN$ is chosen in such a way that
$(C_r \cdot \bN)_\al = 2 N$, for all $\al$) by
$\ch_q \bV_{\bmu}(\bzeta)$. Using the results and notation of section
\ref{awt}, we find
\begin{equation}
\label{chfinf}
\begin{split}
&\ch_q \bV_{\bmu}(\bzeta) =
\sum_{\widetilde{\bm}\in \Z^r}
q^{\frac{1}{2k} \widetilde{\bm}^T C_r \widetilde{\bm}
-\frac{1}{k} \bn^T \cdot \widetilde{\bm}}
\ e^{\bomega^T\cdot \bn-\bomega^T C_{r} \widetilde{\bm}}
 \times \\ & \times  
\frac{1}{(q)_\infty^r}
\sideset{}{'}\sum_{\vbm\in\Z^{r\times(k-1)}_{\geq 0}}
\frac{q^{\frac{1}{2} \vbm^T (C_r \otimes C_{k-1}^{-1}) \vbm
- \vbn^T (\id \otimes C^{-1}_{k-1}) \vbm}}
{\prod_{\al=1}^r\prod_{a<k}(q)_{m_a^{(\al)}}} \ ,
\end{split}
\end{equation}
where the prime denotes the constraint
$\overline{m}^{(\alpha)} = \sum_{a=1}^{k-1} a m^{(\alpha)}_a =  
\widetilde{m}^{(\al)} \mod k$.  As in the case of the fusion of the
principal spaces, the second line of eq. \eqref{chfinf} leads to an
expression for the string functions, in this case associated to
general modules of $\widehat{\mathfrak{sl}}_{r+1}$. However, we can
make the character simpler in appearance by reintroducing $\ma_k$ in
favor of $m^{(\alpha)}$. This gives
\begin{equation}
\label{chfinf2}
\ch_q \bV_{\bmu}(\bzeta) =
\frac{1}{(q)_\infty^r}
\sum_{\substack{\vbm\\\ma_k\in\Z,\ma_{a<k}\in\Z_{\geq 0}}}
\frac{q^{\frac{1}{2}
\vbm^T (C_r \otimes A)\vbm  - \vbn^T (\id \otimes A) \vbm}}
{\prod_{\alpha=1}^{r} \prod_{a=1}^{k-1} (q)_{\ma_a}}
e^{\bomega^T \cdot \bn -\bomega^T C_r \bm}
\ .
\end{equation}

This character decomposes into characters of the integrable modules
in the following way
\begin{equation}\label{fusionint}
\ch_q \bV_{\bmu}(\bzeta) = \sum_{\lambda\leq\nu(\bmu)}
\cK_{\lambda,\bmu} (q^{-1}) \
\ch_q V_{\lambda} \ ,
\end{equation}
where the sum is over dominant weights of $\sl_{r+1}$.


We can now invert the relation \eqref{fusionint}, 
to obtain the character of a general integrable
highest-weight module of $\widehat{\mathfrak{sl}}_{r+1}$.

\begin{thm}
\label{chvl}
The character $\ch_q V_\lambda$ of any integrable, level-$k$
$\widehat{\mathfrak{sl}}_{r+1}$ module with highest weight
$\lambda$ is given by
\begin{equation}
\label{vdecomp}
\ch_q V_{\lambda} = \sum_{\bmu}
(\KK^{-1} (q^{-1}))_{\nu(\bmu),\lambda} \
\ch_q  \bV_\bmu(\bzeta) \ ,
\end{equation}
where $\ch_q  \bV_\bmu(\bzeta)$ is given by equation
\eqref{chfinf2} and the elements of the invertible matrix $\KK$
are given by $(\KK (q))_{\lambda,\nu(\bmu)} =  \cK_{\lambda,\bmu}(q)$,
where $\cK_{\lambda,\bmu}(q)$ is given by equation \eqref{Kostka}.
The finite sum is over sequences of rectangular partitions
of the form $\bmu = (n^{(1)}\omega_1,\ldots,n^{(r)}\omega_r)$, such
that $\nu(\bmu)\leq\lambda$ in the sense of Lemma \ref{kostkaprops}.
\end{thm}

We note some features of this formula. It is a finite sum, with
coefficients in $\Z[q^{-1}]$. Therefore, not only is the 
positivity of the coefficients of $q^n$ not manifest from this
formula, neither is the fact that the character is in fact a series
in positive powers of $q$ only.

\subsection{Some examples}
Let us consider some explicit examples of the matrices of generalized Kostka
polynomials, and, as a result, some character formul\ae\ for
non-rectangular representations.
We will do this for $\widehat{\sl}_3$ in full generality,
and for $\widehat{\sl}_4$ at fixed level. 

\subsubsection{The case $\widehat{\sl}_3$}

In this case, it is very easy to write down the elements of the
matrix $\KK_{\lambda;\nu(\bmu)}$. For a given partition $\lambda$, let
$l_i=\lambda_{i+1}-\lambda_i$. Using this notation, we have the
following result
\begin{equation}
\KK_{(l_1,l_2);(l_1-i,l_2-j)} = \delta_{i,j} q^i \ , 
\end{equation}
where we have the constraints $0 \leq i,j \leq \min(l_1,l_2)$. 
The non-zero elements of $\KK^{-1}$ are also easily obtained
\begin{align}
\KK^{-1}_{(l_1,l_2);(l_1,l_2)} (q)&= 1 \\
\KK^{-1}_{(l_1,l_2);(l_1-1,l_2-1)} (q)&= -q \ , \qquad l_1,l_2 > 0 
\end{align}
while all the other elements are zero. For the characters of arbitrary
$\widehat{\sl}_3$ representations, this implies for non-rectangular
representations (i.e.\ $l_1,l_2 > 0$)
\begin{equation}
\ch_q V_{(l_1,l_2)} = \ch_q \bV_{(l_1,l_2)}(\boldsymbol{\zeta})
- \frac{1}{q} \ch_q\ \bV_{(l_1-1,l_2-1)}(\boldsymbol{\zeta}) \ ,
\end{equation}
where $\ch_q \bV_{\bmu}(\boldsymbol{\zeta})$ is given by equation
\eqref{chfinf} or \eqref{chfinf2}.

\subsubsection{An $\widehat{\sl}_4$ example}

We give an explicit example for the matrix $\KK$ for representations
of $\widehat{\sl}_4$, with level $k\leq 4$. In addition, we will
restrict ourselves to representations with
$\sum_{i=1}^{3} i l_i = 0 \mmod 4$ (see section \ref{kostkastructure}).
There are $10$ representations of
this kind, and we will use the ordering
\begin{equation*}
(0,0,0);(1,0,1),(0,2,0);(2,1,0),(0,1,2);(4,0,0),(2,0,2),(1,2,1),(0,4,0),(0,0,4)
\end{equation*}
With this ordering, we obtain the following Kostka matrix
\begin{equation}
\label{kexample}
\KK (q)=
\left(
\begin{array}{c|cc|cc|ccccc}
1& q&0& 0&0& 0&q^2&0&0&0\\
\hline
0& 1&0& q&q& 0&q&q^2&0&0\\
0& 0&1& 0&0& 0&0&q + q^2&0&0\\
\hline
0& 0&0& 1&0& 0&0&q&0&0\\
0& 0&0& 0&1& 0&0&q&0&0\\
\hline
0& 0&0& 0&0& 1&0&0&0&0\\
0& 0&0& 0&0& 0&1&0&0&0\\
0& 0&0& 0&0& 0&0&1&0&0\\
0& 0&0& 0&0& 0&0&0&1&0\\
0& 0&0& 0&0& 0&0&0&0&1\\
\end{array} \right)
\end{equation}
The inverse is
\begin{equation}
\label{kinvexample}
\KK^{-1}(q) =
\left(
\begin{array}{c|cc|cc|ccccc}
1& -q&0& q^2&q^2& 0&0&-q^3&0&0\\
\hline
0& 1&0& -q&-q& 0&-q&q^2&0&0\\
0& 0&1& 0&0& 0&0&-q - q^2&0&0\\
\hline
0& 0&0& 1&0& 0&0&-q&0&0\\ 
0& 0&0& 0&1& 0&0&-q&0&0\\
\hline
0& 0&0& 0&0& 1&0&0&0&0\\
0& 0&0& 0&0& 0&1&0&0&0\\
0& 0&0& 0&0& 0&0&1&0&0\\
0& 0&0& 0&0& 0&0&0&1&0\\
0& 0&0& 0&0& 0&0&0&0&1\\
\end{array}\right)
\end{equation}

Note that the inverse Kostka matrix has off-diagonal elements with
both signs. As an example, we find that (by making use of equation
\eqref{vdecomp})
\begin{equation}
\begin{split}
\ch_q V_{(1,2,1)} =&\ \ch_q \bV_{(1,2,1)}(\boldsymbol{\zeta})  
- \frac{1}{q} \ch_q \bV_{(2,1,0)}(\boldsymbol{\zeta})
- \frac{1}{q} \ch_q \bV_{(0,1,2)}(\boldsymbol{\zeta}) \\ &
- \bigl(\frac{1}{q}+\frac{1}{q^2}\bigr)
  \ch_q \bV_{(0,2,0)}(\boldsymbol{\zeta})
+ \frac{1}{q^2} \ch_q \bV_{(1,0,1)}(\boldsymbol{\zeta})
- \frac{1}{q^3} \ch_q \bV_{(0,0,0)}(\boldsymbol{\zeta}) \ ,
\end{split}
\end{equation}
with $\ch_q \bV_{\bmu}(\boldsymbol{\zeta})$ given by equation \eqref{chfinf}.

\section{Conclusion}
The main purpose of this paper was to find explicit fermionic
character formul\ae\ for arbitrary integrable highest-weight modules of
$\widehat{\sl}_{r+1}$, using a generalization of the methods of Feigin
and Stoyanovski\u{\i} \cite{FS}. Because the functional realization of the
dual space for non-rectangular highest weights is too complex for
computation of a fermionic character (see Section 3.3.2), we did not
compute purely fermionic characters, which would have the nice feature that
they are manifestly power series in $q$, with non-negative coefficients.
Instead, we found explicit character formul\ae\ as a finite
sum of fermionic characters with coefficients in $\Z[q^{-1}]$.

To obtain these explicit characters, we used the following strategy:
we computed the fermionic character formula for the (non
level-restricted) fusion product of $N$ integrable modules with
rectangular highest weights $\mu_p = a_p \omega_{\alpha_p}$, equation
\eqref{chfinf2}, and of the space of conformal blocks associated with
this fusion product, the generalized Kostka polynomial of Theorem
\ref{generalizedkostka}.  

We thus provided a proof of the conjecture of Feigin and Loktev
\cite{FL}, concerning the relation between their graded tensor product
and the generalized Kostka polynomials \cite{SW,KS} in this case. It
is also a direct proof of the independence of the dimension of the
FL-fusion product of the evaluation parameters (the points $\zeta_p$),
since the associated graded space whose character we computed
corresponds to the limit $\zeta_p\to0$ for all $p$.

We then used the characters for the special case of these fusion
products, together with the relation \eqref{fusionint}, to obtain a
formula for the characters of integrable modules of
$\widehat{\sl}_{r+1}$ of arbitrary (non-rectangular) highest weight,
in terms of the inverse matrix of certain generalized Kostka
polynomials, see Theorem \ref{chvl}.

The generalization of the discussion in this paper to other simple Lie
algebras requires us to consider the so-called Kirillov-Reshetikhin
modules (or rather, their limit to loop algebra case, as
KR-modules were originally defined for Yangians).
These take the place of irreducible $\g$-modules with rectangular highest
weights but as $\g$-modules, they are not necessarily irreducible. We
will explain this generalization in an upcoming publication.

\def\cprime{$'$} \def\cprime{$'$} \def\cprime{$'$}

\end{document}